\documentclass[11pt]{amsart}

\usepackage{amsmath,amssymb,amscd,amsthm,amsxtra}
\usepackage[dvips]{graphics,epsfig}
\usepackage[all]{xy}

\allowdisplaybreaks[2]

\newtheorem{theorem}{Theorem} [section]
\newtheorem{maintheorem}{Theorem}
\newtheorem{lemma}[theorem]{Lemma}
\newtheorem{proposition}[theorem]{Proposition}
\newtheorem{remark}[theorem]{Remark} 

\DeclareMathOperator*{\intt}{\int}

\DeclareMathOperator{\MAX}{MAX}

\newcommand{\noi}{\noindent}

\newcommand{\R}{\mathbb{R}}

\newcommand{\T}{\mathbb{T}}

\newcommand{\al}{\alpha}
\newcommand{\dl}{\delta}
\newcommand{\Dl}{\Delta}
\newcommand{\eps}{\varepsilon}

\newcommand{\g}{\gamma}
\newcommand{\G}{\Gamma}
\newcommand{\ld}{\lambda}

\newcommand{\s}{\sigma}
\newcommand{\ft}{\widehat}
\newcommand{\wt}{\widetilde}
\newcommand{\cj}{\overline}
\newcommand{\dx}{\partial_x}
\newcommand{\dt}{\partial_t}

\newcommand{\I}{\hspace{0.5mm}\text{I}\hspace{0.5mm}}
\newcommand{\II}{\text{I \hspace{-2.8mm} I} }

\newcommand{\jb}[1]
{\langle #1 \rangle}

\numberwithin{equation}{section}
\numberwithin{theorem}{section}

\begin{document}

\title
[ White noise for KdV and mKdV]
{\bf White noise for KdV and mKdV on the circle}

\author{Tadahiro Oh}

\address{Tadahiro Oh\\
Department of Mathematics\\
University of Toronto\\
40 St. George St, Rm 6290,
Toronto, ON M5S 2E4, Canada}

\email{oh@math.toronto.edu}


\subjclass[2000]{ 35Q53}

\keywords{KdV; mKdV; cubic NLS; white noise; invariant measure}

\begin{abstract}
We survey different approaches to study the invariance of the white noise
for the periodic KdV.
We mainly discuss the following two methods.
First, we discuss the PDE method, following  Bourgain \cite{BO4},
in a general framework.  
Then, we show how it can be applied to the low regularity setting of the white noise for KdV
by introducing the Besov-type space $\ft{b}^s_{p, \infty}$, $sp< -1$.
Secondly, we describe the probabilistic method
by Quastel, Valk\'o, and the author \cite{OQV}.
We also use this probabilistic approach  to study the white noise for mKdV.
\end{abstract}

\maketitle


\section{Introduction}

\subsection{KdV, white noise, and Gibbs measure}
In this paper, we  consider the periodic Korteweg-de Vries (KdV) equation:
\begin{equation} \label{KDV}
\begin{cases}
u_t + u_{xxx} +  u u_x  = 0 \\ 
u \big|_{t = 0} = u_0,
\end{cases}
\end{equation}

\noi 
where $u$ is a real-valued function on $\T\times\R$ with $\T = [0, 2\pi)$ and 
the mean of $u_0$ is zero.
By the conservation of the mean, it follows that the solution $u(t)$ of \eqref{KDV} has the spatial mean 0
for all $t \in\mathbb{R}$ as long as it exists.
In the following, we assume that the spatial mean $\ft{u}(0, t)$ is zero for all $t \in \R$.

Our main goal is to show that the mean 0 (Gaussian) white noise  $\mu$ on $\T$
is invariant under the flow of \eqref{KDV}.
Recall that the mean 0  white noise  $\mu$ on $\T$
is the probability measure on distributions $u$ on $\T$ with $\int_\T u dx = 0$
such that we have
\begin{equation}\label{white0}
\int e^{i \jb{ f, u} }d \mu(u) = e^{-\frac{1}{2}\|f\|_{L^2}^2}
\end{equation}

\noi
for all smooth  mean 0 function $f$ on $\T$.
Note that $\jb{ \cdot, \cdot}$ here denotes the $\mathcal{S}-\mathcal{S}'$ duality.
As we see later, we can formally write such $\mu$ as 
\begin{equation}\label{white}
d \mu = Z^{-1}\exp(- \tfrac{1}{2} \int u^2 dx) \prod_{x\in \T} d u(x), \ u \text{ mean } 0 . 
\end{equation}

\noi
From now on, we assume that the spatial mean is always zero,
and hence we may drop the prefix ``mean zero''.

Before providing the precise meaning of \eqref{white}, 
let us discuss the motivation for studying this problem.
 Given a Hamiltonian flow on $\mathbb{R}^{2n}$:
\begin{equation} \label{HR2}
\begin{cases}
\dot{p}_i = \frac{\partial H}{\partial q_i} \\
\dot{q}_i = - \frac{\partial H}{\partial p_i} 
\end{cases}
\end{equation}

\noi
 with Hamiltonian $ H (p, q)= H(p_1, \cdots, p_n, q_1, \cdots, q_n)$,
Liouville's theorem states that the Lebesgue measure on $\mathbb{R}^{2n}$ is invariant under the flow.
Then, it follows from the conservation of the Hamiltonian $H$
that  the Gibbs measures
$e^{-\beta H(p, q)} \prod_{i = 1}^{n} dp_i dq_i$ are invariant under the flow of \eqref{HR2},
where $\beta> 0$ is the reciprocal temperature.

In the context of nonlinear Schr\"odinger equations (NLS), Lebowitz-Rose-Speer \cite{LRS} considered the Gibbs measure of the form 
\begin{equation} \label{Gibbs}
d \mu = \exp (-\beta H(u)) \prod_{x\in \T} d u(x),
\end{equation}

\noi
where 
$H(u)$ is the Hamiltonian given by $H(u) = \frac{1}{2} \int |u_x|^2 \pm \frac{1}{p} \int |u|^p dx$.
They showed that such Gibbs measure $\mu$ is  
a well-defined probability measure on $H^{\frac{1}{2}-}(\T) := \bigcap_{s<\frac{1}{2}} H^s(\T)$.
(In the focusing case (with  $-$), 
the result only holds for $p < 6$ 
with the $L^2$-cutoff $\chi_{\{\|u\|_{L^2} \leq B \}}$ for any $B>0$, 
and for $ p = 6$ with sufficiently small $B$.)
Using the Fourier analytic approach, Bourgain \cite{BO4} continued this study and 
proved the invariance of the Gibbs measure $\mu$ under the flow of NLS
and global well-posedness almost surely on the statistical ensemble.
He also established the invariance of the Gibbs measures
for KdV, mKdV \cite{BO4}, the Zakharov system on $\T$ \cite{BO5}, 
defocusing cubic NLS in $\T^2$ and $\T^3$ \cite{BO7, BO8},
defocusing cubic NLS on $\R$ \cite{BO9}.

There are many results in this direction:
Friedlander \cite{FR}, Zhidkov \cite{Z1, Z2}, 
McKean-Vaninsky \cite{MV1, MV2}, McKean \cite{MK1, MK2} for NLS and 
nonlinear wave equations (NLW) on $\T$ (and on $\R$ \cite{MV3}.)
Note that some of them employ non-Fourier analytic technique,
and they are rather probabilistic.
There are more recent results based on Bourgain's approach:
Tzvetkov \cite{TZ1, TZ2} for subquintic radial NLS on the unit disc in $\R^2$, 
Burq-Tzvetkov \cite{BT1, BT3} for subquartic NLW on the unit ball in $\mathbb{R}^3$ with radial symmetry and
the Dirichlet boundary condition,
Oh \cite{OH3} for the coupled KdV systems under certain Diophantine conditions 
and \cite{OHSBO} for the Schr\"odinger-Benjamin-Ono system.

This study was partially motivated to answer the question posed 
by V. Zakharov during the Sixth I. G. Petrovskii memorial meeting of the Moscow Mathematical
Society in 1983. c.f. \cite{FR}.
``Numerical experiments demonstrated [that the 1-d periodic cubic NLW]
possesses the ``returning" property, i.e. solutions appear to be very close to the initial state
$\cdots$, after some time of rather chaotic evolution. 
The problem is to explain  this phenomenon.'' 
Also, see the related Fermi-Pasta-Ulam phenomenon \cite{Palais}.

Invariant Gibbs measures $\mu$ for the nonlinear PDEs can be regarded
as invariant measures  for infinite dimensional dynamical systems,
and it follows from  Poincar\'e recurrence theorem 
that almost all the points of the phase space  
are stable according to Poisson. See Zhidkov \cite{Z}.
Note that this recurrence property holds only in the support of the Gibbs measure,
i.e. not for smooth functions, except for the integrable PDEs as mentioned below.

Lastly, note that if $F(p, q)$ is any function that is conserved under the flow of \eqref{HR2}, 
then 
the measure $d \mu_F = 
e^{-\beta F(p, q)} \prod_{i = 1}^{n} dp_i dq_i$ is invariant.
For KdV and cubic NLS, Zhidkov \cite{Z3, Z4} constructed infinite sequences
of the invariant measures on smoother phase spaces corresponding 
to the higher order conserved quantities of these integrable PDEs.
Now, recall that KdV \eqref{KDV} preserves the $L^2$ norm.
Hence, it seems natural, at least at a heuristic level, 
to expect the invariance of the white noise \eqref{white}.
The difficulty here is the low regularity of the phase space as we see in the next subsection.

\subsection{Gaussian measures in Hilbert spaces}
In this subsection,  
we briefly go over the basic theory of Gaussian measures in Hilbert spaces
to provide the precise meaning of \eqref{white}.
See  Zhidkov \cite{Z} for details. 

First, recall  (centered) Gaussian measures in $\mathbb{R}^n$.
Let $n \in \mathbb{N}$ and $B$ be a symmetric positive $n \times n$ matrix
with real entries.
The  Borel measure $\mu$ in $\mathbb{R}^n$ with the density
\[ d \mu(x) = \frac{1}{\sqrt{(2\pi)^n \det (B )}} \exp \big( -\tfrac{1}{2} \langle B^{-1} x, x \rangle_{\mathbb{R}^n} \big)\]

\noindent
is called a (nondegenerate centered) Gaussian measure in $\mathbb{R}^n$.
Note that $\mu(\mathbb{R}^n) = 1$.

Now, we consider an analogous definition of infinite dimensional (centered) Gaussian measures.
Let $H$ be a real separable Hilbert space and $B: H \to H$ be a linear positive self-adjoint operator 
(generally not bounded) with eigenvalues $\{\ld_n\}_{n\in \mathbb{N}}$
and the corresponding eigenvectors $\{e_n\}_{n\in\mathbb{N}}$  forming an orthonormal basis of $H$.
We call a set $M \subset H$  cylindrical if there exists an integer $n\geq 1$ and a Borel set $F \subset \mathbb{R}^n$
such that
\begin{equation} \label{CYLINDER}
 M = \big\{ x \in H : ( \jb{ x, e_1}_H, \cdots, \jb{ x, e_n}_H ) \in F \big\}. 
\end{equation}

\noindent
For a fixed operator $B$ as above, we denote by $\mathcal{A}$ the set of all cylindrical subsets of $H$.
One can easily verify that $\mathcal{A}$ is a field.
Then, the centered Gaussian measure in $H$ with the correlation operator $B$ is defined as 
the additive (but not countably additive in general) measure $\mu$ defined on the field $\mathcal{A}$
via
\begin{equation} \label{CGAUSSIAN}
 \mu(M) = (2\pi)^{-\frac{n}{2}} \prod_{j = 1}^n \ld_j^{-\frac{1}{2}} \int_F e^{-\frac{1}{2}\sum_{j = 1}^n \ld_j^{-1} x_j^2 }d x_1 \cdots dx_n,
\text{ for }M \in \mathcal{A} \text{ as in \eqref{CYLINDER}. }
\end{equation}

\noi
The following proposition tells us when  this Gaussian measure $\mu$ is countably additive.

\begin{proposition} \label{COUNTABLEADD}
The  Gaussian measure $\mu$ defined in \eqref{CGAUSSIAN} is countably additive
on the field $\mathcal{A}$ if and only if $B$ is an operator of trace class, 
i.e. $\sum_{n = 1}^\infty \ld_n < \infty$.
If the latter holds, then
the minimal $\s$-field $\mathcal{M}$ containing the field $\mathcal{A}$ of all cylindrical sets is the Borel $\s$-field on $H$.
\end{proposition}

Consider a sequence of the finite dimensional Gaussian measures $\{\mu_n\}_{n\in\mathbb{N}}$ as follows.
For fixed $n \in \mathbb{N}$, let $\mathcal{M}_n$ be the set of all cylindrical sets in $H$ of the form \eqref{CYLINDER} with this fixed $n$
and arbitrary Borel sets $F\subset \mathbb{R}^n$.
Clearly, $\mathcal{M}_n$ is a $\s$-field, and setting 
\[ \mu_n(M) = (2\pi)^{-\frac{n}{2}} \prod_{j = 1}^n \ld_j^{-\frac{1}{2}} \int_F e^{-\frac{1}{2}\sum_{j = 1}^n \ld_j^{-1} x_j^2 }d x_1 \cdots dx_n\]

\noindent
for $M \in \mathcal{M}_n$, we obtain a countably additive measure $\mu_n$ defined on $\mathcal{M}_n$.
Then, one can show that each measure $\mu_n$ can be naturally extended onto the whole Borel $\s$-field $\mathcal{M}$ of $H$
by $ \mu_n(A) := \mu_n(A \cap \text{span}\{e_1, \cdots, e_n\})$
for $A \in \mathcal{M}$.
Then, we have

\begin{proposition} \label{PROP:Zhidkov2}
Let $\mu$ in  \eqref{CGAUSSIAN} be countably additive.
Then,  $\{\mu_n\}_{n\in \mathbb{N}}$ constructed above converges weakly to $\mu$ as $n \to \infty$.

\end{proposition}

Now,  we construct the mean 0 white noise.
Let $u = \sum_{n} \ft{u}_n e^{inx}$ be a real-valued function on $\mathbb{T}$ with mean 0. 
i.e. we have $\ft{u}_0 = 0$ and $\ft{u}_{-n} = \cj{\ft{u}_n}$.
First, define $\mu_N$ on $\mathbb{C}^{N} \cong \mathbb{R}^{2N}$
with the  density
\begin{equation}\label{WhiteN}
 d \mu_N = Z_N^{-1} e^{- \sum_{n = 1}^N  |\ft{u}_n|^2 } \textstyle \prod_{ n = 1 }^N d \ft{u}_n , 
\end{equation}

\noindent
where
$Z_{N} = \int_{\mathbb{C}^{N}} 
e^{- \sum_{ n = 1 }^N  |\ft{u}_n|^2} \prod_{ n = 1 }^N  d \ft{u}_n . $
Note that this measure is the induced probability measure on $\mathbb{C}^{N}$ under the map
\begin{equation} \label{induced}
 \omega \mapsto \{  g_n(\omega)  \}_{n = 1}^N,
\end{equation}

\noindent
where $g_n(\omega)$, $n = 1, \cdots, N$,  are independent standard complex Gaussian random variables.
Next, define the white noise $\mu$ by
\begin{equation} \label{White}
 d \mu = Z^{-1} e^{- \sum_{  n \geq 1 }  |\ft{u}_n|^2} \textstyle \prod_{  n \geq 1 } d  \ft{u}_n , 
\end{equation}

\noindent
where
$Z = \int e^{- \sum_{   n \geq 1 }  |\ft{u}_n|^2}  \prod_{  n \geq 1 } d \ft{u}_n. $
Then, in the above correspondence, we have 
$u = \sum_{n \ne 0} g_n e^{inx}$, 
where $\{g_n(\omega)\}_{n \geq 1}$ are independent standard complex Gaussian random variables
and $g_{-n} = \cj{g_n}$.

Let $\dot{H}^s_0$ be the homogeneous Sobolev space restricted to the {\it real-valued} mean 0 elements.
Let $\jb{\cdot, \cdot}_{\dot{H}^s_0}$ denote the inner product in $\dot{H}^s_0$.
i.e. $ \big\langle \sum c_n e^{inx}, \sum d_n e^{inx} \big\rangle_{\dot{H}_0^s} 
= \sum_{n \ne 0} |n|^{2s} c_n \cj{d_n} $.
Let $B_s = \sqrt{-\Dl}\vphantom{|}^{2s}$.
Then, the weighted exponentials $\{|n|^{-s} e^{inx}\}_{n\ne 0}$ are the eigenvectors of $B_s$ with the eigenvalue $|n|^{2s}$,
forming an orthonormal basis of $\dot{H}^s_0$.
Note that
\[ -\tfrac{1}{2}\jb{B^{-1} u, u}_{\dot{H_0^s}} 
= -\tfrac{1}{2}\Big\langle \sum_{n \ne 0} |n|^{-2s} \ft{u}_n e^{inx}, 
\sum_{n \ne 0} \ft{u}_n e^{inx} \Big\rangle_{\dot{H}_0^s} 
= -\sum_{n \geq 1}  |\ft{u}_n|^2.\] 

\noindent
The right hand side is exactly the expression appearing in the exponent in \eqref{White}.
It follows from Proposition \ref{COUNTABLEADD} that $\mu$ is  countably additive
if and only if $B$ is  of trace class, i.e.
$ \sum_{n \ne 0} |n|^{2s} < \infty$.
Hence, $\bigcap_{s < -\frac{1}{2}} H^s$ is a natural space to work on.

\begin{remark} \rm
In view of \eqref{induced} with $N = \infty$, 
we see that $u$ in the support of the white noise \eqref{white}
has the representation $u = \sum_{n \ne 0} g_n(\omega) e^{inx}$.
Then, for a smooth  mean 0 function $f$ on $\T$, we have
\begin{equation}
\int e^{i \jb{ f, u} }d \mu(u) = \prod_{n\ne0} \int e^{i \ft{f}_n \cj{\ft{g}_n}} d g_n
= e^{-\sum_{n\geq 1 } |\ft{f}_n|^2}
= e^{-\frac{1}{2}\|f\|_{L^2}^2}.
\end{equation}

\noi
Hence, \eqref{white0} is satisfied.

Moreover, we can regard $u \in \text{supp}(\mu)$ as the Gaussian
randomization of the Dirac delta $\delta_0(x)$ on the Fourier coefficients.
Recall that $\delta_0(x)$ is in $H^s(\T)$ for $s < -\frac{1}{2}$
but not in $H^{-\frac{1}{2}}(\T)$.
It is also known \cite{BT2} that the Gaussian randomization of the Fourier coefficients
does {\it not} give any smoothing (in terms of the Sobolev regularity) a.s.
This also shows that 
$\text{supp}(\mu) \subset \bigcap_{s < -\frac{1}{2}} H^s \setminus H^{-\frac{1}{2}}$.

\end{remark}

\section{Main results}

In this section, we state several different methods for proving the invariance of the white noise.
They are arranged in the chronological order,
and Methods 1 and 2 are described more in details in the following sections.

\subsection{Method 0: Complete integrability approach}

The first result for the invariance of the white noise for \eqref{KDV}
is due to Quastel-Valk\'o \cite{QV}.
This exploits the bi-Hamiltonian structure of the KdV:
\begin{equation} \label{KDV2}
u_t + u_{xxx} - 6 u u_x = 0.
\end{equation}

\noi
Recall that \eqref{KDV2} can be written as $u_t = J_i \frac{d H_i}{du}$, $i = 1, 2$,
where $J_1 = \dx$ and $H_1 =  \int \frac{1}{2} u_x^2 + u^3$, i.e. the usual Hamiltonian structure, 
and $J_2 = \dx^3 + 4 u \dx + 2 \dx u$ and $H_2 = \int u^2$.

Their argument uses the correspondence between $(J_2, H_2)$ for KdV 
and the usual Hamiltonian structure of mKdV:
\begin{equation} \label{mKDV}
\begin{cases}
u_t + u_{xxx} \pm  u^2 u_x  = 0 \\ 
u \big|_{t = 0} = u_0,
\end{cases}
\end{equation}

\noi
More precisely, their argument combines the following results:
\begin{enumerate}
\item[(i)]
Cambronero-McKean \cite{CM}:
the (corrected) Miura transform maps
the usual Gibbs measure of the form \eqref{Gibbs} for mKdV \eqref{mKDV}
(with the $-$ sign)
to  the mean 0 white noise \eqref{white} for KdV.

\item[(ii)]
Bourgain \cite{BO4}: invariance of the Gibbs measure for mKdV.

\item[(iii)]
Kappeler-Topalov \cite{KT}:  global well-posedness (GWP) of \eqref{KDV} in $H^{-1}(\T)$ 
via the inverse spectral method. 

\end{enumerate}

\noi
(i) and (ii) imply that the white noise \eqref{white} for KdV
is invariant if the flow is well-defined in its support,
and (iii) guarantees such well-posedness.
Note that this method heavily depends on the complete integrability of \eqref{KDV}
and is not applicable to the general non-integrable variants of KdV, 
including the coupled KdV system considered in \cite{OH3}.

\subsection{Method 1: Bourgain's PDE approach}

First, note that the invariance of the white noise follows
once we show that 
\eqref{KDV} is almost surely globally well-posed with $u_0 = \sum_{n\ne0} g_n(\omega)e^{inx}$
and that $u(t)$ has the same distribution for all $t\in \R$.

In \cite{BO4}, Bourgain proved the invariance of the Gibbs measures for NLS.  
In dealing with super-cubic nonlinearity,
(where only local well-posedness (LWP)  was available),
he used a  probabilistic argument and the approximating finite dimensional ODEs
(with the invariant finite dimensional Gibbs measures)
to extend the local solutions to global ones almost surely on the statistical ensembles.
Then, he proved the invariance of the Gibbs measures.
Note that it was crucial that LWP was obtained with a ``good" estimate
on the solutions for his argument to obtain the uniform convergence of the solutions 
of the finite dimensional ODEs to those of the full PDE.
e.g. Lemma 41 in \cite{BO6}.
For the details of the argument, 
see Bourgain \cite[Lec.4]{BO6}, Burq-Tzvetkov \cite[Sec.6-7]{BT1}, Oh \cite[Sec.6]{OH3}, and 
Tzvetkov \cite[Sec.8-10]{TZ1}, \cite[Sec.7-9]{TZ2}.

Hence, the main difficulty in  this approach is to establish LWP of \eqref{KDV}
in the support of the white noise $\mu$.
Then, we can establish the invariance by following Bourgain's argument.

Now, we briefly review recent well-posedness results of the periodic KdV \eqref{KDV}.
In \cite{BO1}, Bourgain  introduced a new weighted space-time Sobolev space $X^{s, b}$
whose norm is given by
\begin{equation} \label{Xsb}
\| u \|_{X^{s, b}(\mathbb{T} \times \mathbb{R})} = \| \jb{n}^s \jb{\tau - n^3}^b 
\ft{u}(n, \tau) \|_{L^2_{n, \tau}(\mathbb{Z} \times \R)},
\end{equation}

\noindent
where $\jb{ \: \cdot \:} = 1 + |  \cdot  | $. 
He proved local well-posedness of \eqref{KDV} in $L^2(\mathbb{T})$
via the fixed point argument, 
immediately yielding global well-posedness in $L^2(\mathbb{T})$
thanks to the conservation of the $L^2$ norm.
Kenig-Ponce-Vega \cite{KPV4} (also see \cite{CKSTT4}) improved Bourgain's result 
and established local well-posedness in $H^{-\frac{1}{2}}(\T)$
by establishing the bilinear estimate 
\begin{equation} \label{KPVbilinear}
\| \dx(uv) \|_{X^{s, -\frac{1}{2}}} \lesssim \| u \|_{X^{s, \frac{1}{2}}} \| v \|_{X^{s, \frac{1}{2}}}, 
\end{equation}

\noi
for $s \geq -\frac{1}{2}$ under the mean 0 assumption on $u$ and $v$.
Colliander-Keel-Staffilani-Takaoka-Tao \cite{CKSTT4} proved 
the corresponding global well-posedness result via the $I$-method. 

There are also results on \eqref{KDV} which exploit its complete integrability.
In \cite{BO3}, Bourgain proved global well-posedness of \eqref{KDV}
in the class $M(\T)$ of measures $\mu$, assuming that 
its total variation $\|\mu\|$ is sufficiently small.
His proof is based on the trilinear estimate on the second iteration of 
the integral formulation of \eqref{KDV}, 
assuming an a priori uniform bound on the Fourier coefficients of the solution $u$ of the form 
\begin{equation} \label{BOO}
\sup_{n\in \mathbb{Z}} |\ft{u}(n, t)| < C
\end{equation}

\noi
for all $t\in \R$. 
Then, he established \eqref{BOO} using the complete integrability.
More recently, Kappeler-Topalov \cite{KT} proved global well-posedness of the KdV in $H^{-1}(\T)$
via the inverse spectral method as already mentioned in Subsection 2.1. 

There are also results on necessary conditions on the  regularity 
with respect to smoothness or uniform continuity 
of the solution map $: u_0 \in H^s (\mathbb{T}) \to u(t) \in H^s(\mathbb{T})$.
Bourgain \cite{BO3} showed that if the solution map is $C^3$, 
then $s \geq -\frac{1}{2}$.  
Christ-Colliander-Tao \cite{CCT}
proved that if the solution map is uniformly continuous, 
then $s \geq -\frac{1}{2}$.
(Also, see Kenig-Ponce-Vega \cite{KPV5}.) 
These results, in particular, imply that
we can not hope to have a local-in-time solution of \eqref{KDV} via the fixed point argument in $H^s$, $s < -\frac{1}{2}$.

Recall that
the white noise $\mu$ defined in \eqref{white} is supported on $\cap_{s < -\frac{1}{2}}H^s \setminus H^{-\frac{1}{2}}$ a.s.
Hence, we can not simply apply the known results to study local well-posedness of \eqref{KDV}
in the support of $\mu$.
Instead, we  prove a local well-posedness in an appropriate Banach space containing the support of the white noise $\mu$.
Define a Besov-type space $\ft{b}^s_{p, \infty}$
 via the norm
\begin{equation} \label{Besov}
\| f\|_{\ft{b}^s_{p, \infty}} 
:= \| \ft{f}\|_{b^s_{p, \infty}} = \sup_j \| \jb{n}^s \ft{f}(n) \|_{L^p_{|n|\sim 2^j}}
= \sup_j \Big( \sum_{|n| \sim 2^j} \jb{n}^{sp} |\ft{f}(n)|^p \Big)^\frac{1}{p}.
\end{equation}

\noi
By Hausdorff-Young's inequality, 
we have $\ft{b}^s_{p, \infty} \supset B^s_{p', \infty}$ for $p >2$,
where $B^s_{p', \infty}$ is the usual Besov space with $p' = \frac{p}{p-1}$.
This space has two important properties:
\begin{enumerate}
\item[(i)]
$\ft{b}^s_{p, \infty}$ contains the support of the white noise for $sp < -1$.
This follows from  the theory of abstract Wiener spaces (c.f. Gross \cite{GROSS}, Kuo \cite{KUO}.)
See  \cite{OH4}.

\item[(ii)]
We can carry out the nonlinear analysis on the second iteration introduced by Bourgain in \cite{BO3}, 
{\it without} assuming the a priori bound \eqref{BOO},
if we take the initial data $u_0 \in \ft{b}^s_{p, \infty}$ for  $s >-\frac{1}{2}$
with $p > 2$.
Then, we construct a solution $u$ as a strong limit of the smooth solutions $u^{(n)}$
 of \eqref{KDV} with smooth $u_0^{(n)}$.
See \cite{OH6}. 
\end{enumerate}

\noi
Hence, we establish LWP in a Banach space containing the support of $\mu$.

\begin{maintheorem} \label{THM:LWP2}
Assume the mean 0 condition on $u_0$.
Let $p = 2+$ and $s = -\frac{1}{2}+ \dl$ with $\frac{p-2}{4p} <\dl < \frac{p-2}{2p}$.
i.e. $s p < -1$.
Then, \eqref{KDV} is locally well-posed in $\ft{b}^s_{p, \infty}$.
\end{maintheorem}

\noi
Although this LWP is not obtained via the fixed point argument, 
the estimates are strong enough to conclude a.s. GWP and the invariance of the white noise,
following Bourgain's argument.
We describe some of the details in Section 3.

\begin{maintheorem} \label{THM:GWP2}
Let $\{ g_n (\omega) \}_{n = 1}^\infty$ be a sequence of independent standard  complex Gaussian
random variables on a probability space $(\Omega, \mathcal{F}, P)$. 
Consider \eqref{KDV} with initial data 
$ u_0 = \sum_{n \ne 0} g_n (\omega) e^{inx},$
where $g_{-n} = \cj{g_n}$.
Then, \eqref{KDV} is globally well-posed almost surely in $\omega \in \Omega.$
Moreover, the mean 0 white noise $\mu$ is invariant under the flow. 
\end{maintheorem}

\begin{remark}\rm
Theorem \ref{THM:LWP2} provides an answer to the question posed by Bourgain 
in \cite[Remark on p.120]{BO3}, 
at least in the local-in-time setting.
i.e. it establishes local well-posedness of \eqref{KDV}
for  a finite Borel measure  $u_0 = \mu \in M(\T)$
with $\|\mu\| < \infty$ {\it without} the complete integrability or the smallness assumption on $\|\mu\|$.
Just note that $\mu \in \ft{b}^s_{p, \infty}$
for $sp \leq -1$ since $\sup_n |\ft{\mu}(n)| <\|\mu\|< \infty.$
Hence, it can be used to  study the Cauchy problem on $M(\T)$ for non-integrable KdV-variants.
\end{remark}

\begin{remark}
Let $\mathcal{F} L^{s, p}$ be the space of functions on $\mathbb{T}$ defined via the norm
\begin{equation} \label{FLP}
\|f\|_{\mathcal{F} L^{s, p}}=\| \jb{n}^s \ft{f}(n) \|_{L^p_n}.
\end{equation}

\noi
Then, Theorems \ref{THM:LWP2} and \ref{THM:GWP2} can also be established in 
$\mathcal{F} L^{s, p}$
with 
$p = 2+$ and $s = -\frac{1}{2}+ \dl$ with $\frac{p-2}{4p} <\dl < \frac{p-2}{2p}$.
\end{remark}

\subsection{Method 2: Probabilistic approach}

Now, we discuss the probabilistic approach by Oh-Quastel-Valk\'o \cite{OQV}.
In this approach, we consider the interpolation of the Gibbs measure of the form \eqref{Gibbs} and the white noise \eqref{white}.
First, consider the Gaussian measure $\mu_\beta$ given by 
\begin{align} \label{Gauss2}
d \mu_\beta =  Z_\beta^{-1}
e^{-\frac{1}{2} \int u^2 - \frac{\beta}{2} \int u_x^2}
\prod_{x \in \T} du(x), 
\end{align}

\noi
where $u$ is real-valued with mean 0.
This is an interpolation of the Wiener measure and the white noise on $\T$.
In the support of $\mu_\beta$,
$u$ has the representation:
\begin{equation} \label{representation1}
u(x) = \sum_{n \ne 0} \frac{g_n}{\sqrt{1  + \beta n^2}} e^{2\pi i n x}, 
\ g_{-n} = \cj{g_n}.
\end{equation}

\noi
i.e.
for each $\beta > 0$, 
$u$ is a.s. in $H^s$ for $s <\frac{1}{2}$ but not in  $H^{\frac{1}{2}}$.
When $\beta = 0$, 
\eqref{Gauss2} reduces to the usual white noise \eqref{white}
supported on $\bigcap_{s<-\frac{1}{2}} H^s \setminus H^{-\frac{1}{2}}$.

Now, define the interpolation  of the Gibbs measure and the white noise by
\begin{align} \label{Gibbs2}
d \rho_\beta  = d \rho^{(p)}_{ \beta} & := 
\wt{Z}_\beta^{-1} \chi_{\{ \int u^2 \leq K \beta^{-\frac{1}{2}} \}} e^{\beta \int u^p} d \mu_\beta \notag \\
& \hphantom{:}= \ft{Z}_\beta^{-1}\chi_{\{ \int u^2 \leq K \beta^{-\frac{1}{2}} \}}
e^{-\frac{1}{2} \int u^2 + \beta \int u^p - \frac{\beta}{2} \int u_x^2}
\prod_{x \in \T} du(x).
\end{align}

\noi
From \cite{LRS} and \cite{BO4}, 
we see that $e^{\beta \int u^p}$ is integrable with respect to $\mu_\beta$
for each fixed $\beta > 0$ and $p < 6$.
Hence, we can choose appropriate normalizing constants
$\wt{Z}_\beta = \wt{Z}_\beta(p)$ and $\ft{Z}_\beta= \ft{Z}_\beta(p)$ 
so that $\rho_\beta = \rho_\beta^{(p)}$
is a probability measure.

Set $p = 3$  for KdV.
For each $\beta > 0$, $\rho_\beta^{(3)}$ basically behaves like 
the Gibbs measure for KdV and thus it is invariant under the KdV flow.
Moreover, we have the following weak convergence result.

\begin{maintheorem} \label{THM:WC1}
$\rho_\beta^{(3)}$ converges weakly to the white noise $\mu$ as $\beta \to 0$.
\end{maintheorem}

\noi
Hence, 
the white noise $\mu$ is a weak limit of the invariant measures $\rho^{(3)}_\beta$,
and we expect such a measure to be invariant as well.
In this case, we can establish the invariance of the white noise $\mu$
thanks to  the continuity of the KdV flow in $H^{-1}(\T)$ containing the support of $\mu$.
Note that it is enough to have the continuity of the flow in the support.
i.e. Unlike the PDE approach, 
we  do not need any estimate.

The main difficulty of the proof of Theorem \ref{THM:WC1}
lies in establishing the exponential expectation estimate:
\begin{equation} \label{ExpEx1}
\mathbb{E}_{\mu_\beta} \big[\chi_{\{ \int u^2 \leq K \beta^{-\frac{1}{2}} \}} e^{r \beta \int u^3}\big]
= \int \chi_{\{ \int u^2 \leq K \beta^{-\frac{1}{2}} \}} e^{r \beta \int u^3} d \mu_\beta
\leq C(r) <\infty,
\end{equation}

\noi 
uniformly in $\beta > 0$.

Now, let's turn out attention to mKdV, i.e. $p = 4$.
As before, $\rho_\beta^{(4)}$ behaves like 
the Gibbs measure for mKdV for each $\beta > 0$. 
Thus,  Bourgain's result on the invariance of the Gibbs measure for mKdV
implies that $\rho_\beta^{(4)}$ is invariant under the mKdV flow for $\beta > 0$.
Moreover, we can prove the weak convergence of $\rho_\beta^{(4)}$ to the white noise
in this case as well.

\begin{maintheorem} \label{THM:WC2}
$\rho_\beta^{(4)}$ converges weakly to the white noise $\mu$ as $\beta \to 0$.
\end{maintheorem}

\noi
Unfortunately, this does not establish the invariance of the white noise for mKdV,
since the flow of mKdV is not well-defined in the support of the white noise.
(Recall that mKdV is scaling-supercritical in $H^s$ for $s < -\frac{1}{2}$.) 
Theorem \ref{THM:WC2} implies a version of ``formal'' invariance of the white noise in the following sense.

Let $u^\beta_0 $ be a random variable on $\T$ with distribution $\rho_\beta^{(4)}$.
The solution $u^\beta(t, \omega)$ of mKdV with
$u^\beta(0) = u^\beta_0 (\omega)$ exists globally in time,
almost surely in  $\omega \in  \Omega$.
Moreover, $u^\beta(t, \omega)$ has the same distribution $\rho_\beta^{(4)}$ for all $t \in \R$.
By Theorem \ref{THM:WC2},
$u^\beta_0 $ converges weakly to $u_0$, a random variable with the white noise $\mu$ as its distribution.
Also, for each $t \in \R$, $u^\beta (t)$ converges weakly to 
{\it some} random variable $v_t$ with $\mu$ as its distribution.
We would like to say that $v_t = u(t)$, the solution of mKdV 
with $u(0) = u_0 (\omega)$, which would then imply the invariance of the white noise for mKdV. 
However, the flow of mKdV in the support of the white noise
is not known to be well-defined.

\begin{remark} \label{RM:2.3} \rm
Recall that mKdV is scaling-supercritical in $H^s$ for $s < -\frac{1}{2}$,
and the support of the white noise is 
contained in $\bigcap_{s<-\frac{1}{2}} H^s \setminus H^{-\frac{1}{2}}$.
This does {\it not} imply that it is impossible to define the flow on the support of the white noise.
Indeed, we may be able to define the flow of mKdV just on the support of the white noise.
See Bourgain \cite{BO7}
for the case of the 2-$d$ defocusing cubic NLS. 
Recall that the 2-$d$ cubic NLS is $L^2$-critical,
and the Wiener measure on $\T^2$ (the Gaussian part of
the Gibbs measure) is supported below $L^2 (\T^2)$.
Nonetheless, Bourgain constructed a well-defined flow on its support
(after the Wick ordering on the nonlinearity -- a kind of renormalization related to 
the Euclidean quantum field theory),
and established the invariance of the Gibbs measure.
\end{remark}

As in the $p = 3$ case, 
the main difficulty of the proof of Theorem \ref{THM:WC2}
lies in establishing the exponential expectation estimate:
\begin{equation} \label{ExpEx2}
\mathbb{E}_{\mu_\beta} \big[\chi_{\{ \int u^2 \leq K \beta^{-\frac{1}{2}} \}} e^{r \beta \int u^4}\big]
= \int \chi_{\{ \int u^2 \leq K \beta^{-\frac{1}{2}} \}} e^{r \beta \int u^4} d \mu_\beta
\leq C(r) <\infty,
\end{equation}

\noi 
uniformly in $\beta > 0$.
It turns out that \eqref{ExpEx2} is much more delicate than \eqref{ExpEx1}.
We need some probabilistic tools such as
the hypercontractivity of the Ornstein-Uhlenbeck semigroup.
We discuss some of the details in Section 4.
Lastly, we point out that a result similar to Theorem \ref{THM:WC2}
holds for the 1-$d$ cubic NLS, which is also $H^{-\frac{1}{2}}$-critical.
Once again, this result establishes only the formal invariance of the white noise
in the sense  described above.

\section{Method 1: Bourgain's PDE approach}
\subsection{General framework}
In this subsection, we review Bourgain's idea
in a general framework,
and discuss how to prove almost surely GWP and the invariance of a measure from LWP.
Consider a dispersive nonlinear Hamiltonian PDE with a $k$-linear nonlinearity:
\begin{equation} \label{PDE1}
\begin{cases}
u_t = \mathcal{L}u + \mathcal{N}(u)\\
u|_{t = 0} = u_0
\end{cases}
\end{equation}

\noi
where $\mathcal{L}$ is a (spatial) differential operator
whose symbol $P(\xi)$ is  given by a polynomial  with real coefficients on the odd degree terms 
and purely imaginary coefficients on the even degree terms,
and $\mathcal{N}(u) = \mathcal{N}(u, \cdots, u)$ is a $k$-linear nonlinearity, 
possibly with a derivative.
Let $H(u)$ denote the Hamiltonian of \eqref{PDE1}.
Then, \eqref{PDE1} can also be written 
as $u_t = J \, \frac{d H}{d u}$ if $u$ is real-valued,
and as $u_t = J \,  \frac{\partial H}{\partial \cj{u}}$ if $u$ is complex-valued.
For simplicity, we assume that $u$ is real-valued.

Let $\mu$ denote a measure on 
the distributions on $\T$,
whose invariance we'd like to establish.
We assume that $\mu$ is a (weighted) Gaussian measure
given by $d\mu = Z^{-1} e^{-F(u)} \prod_{x \in \T} d u(x)$,
where $F(u)$ is conserved under the flow of \eqref{PDE1}
and the leading term of $F(u)$ is quadratic and nonnegative.

Now, suppose that there exist a Banach space $B$ of distributions on $\T$
and 
a space $X_\dl$ of space-time distributions 
such that we have the following:
\begin{enumerate}
\item[(i)] $X_\dl \subset C([-\dl, \dl]; B)$, 
and $\text{supp}(\mu) \subset B$ 
in the sense that  
$(B, \mu)$ is an abstract Wiener space. See Subsection 3.2.

\vspace{3pt}
\item[(ii)] linear homogeneous estimate:
$\| S(t) u_0\|_{X_\dl} \lesssim \|u_0\|_B$, where $S(t) = e^{t \mathcal{L}}$ 
\vspace{3pt}
\item[(iii)] linear nonhomogeneous estimate:
$ \|\int_0^t S(t - t') F(t') dt'\|_{X_\dl} \lesssim \|F\|_{X'_\dl}$
\vspace{1pt}
\item[(iv)] $k$-linear estimate:
$\|\mathcal{N}(u_1, \cdots, u_k )\|_{X'_\dl} \lesssim \dl^\theta 
\prod_{j = 1}^k \|u_j\|_{X_\dl}$
\end{enumerate}

\noi
for some appropriate auxiliary space $X'_\dl$ and $\theta > 0$.
Then, it is easy to see that 
\eqref{PDE1} is LWP via the fixed point argument 
in the subcritical sense.
i.e. the time of local existence $\dl$ depends on $\|u_0\|_B$,
say $\dl \sim \|u_0\|_B^{-\al}$ for some $\al > 0$.

In addition, suppose that the Dirichlet projections $\mathbb{P}_{ N}$ --
the projection onto the spatial frequencies $\leq N$ --
act boundedly on these spaces, uniformly in $N$.
Consider the finite dimensional approximation to \eqref{PDE1}:
\begin{equation} \label{PDE2}
\begin{cases}
u^N_t =  \mathcal{L}u^N + \mathbb{P}_N\big(\mathcal{N}(u^N)\big)\\
u^N|_{t = 0} = u^N_0
:= \mathbb{P}_{ N} u_0 (x) = \sum_{|n| \leq N} \ft{u}_0(n) e^{inx}.
\end{cases}
\end{equation}

\noi
Then, for $\|u_0\|_B \leq K$, it follows that 
\eqref{PDE2} is LWP on $[-\dl, \dl]$ with 
$\dl \sim K^{-\al}$, independent of $N$.
We need two more assumptions on \eqref{PDE2}.
\begin{enumerate}
\item[(v)] \eqref{PDE2} is Hamiltonian with $H(u^N)$. 
i.e. $u^N_t = J \, \frac{d H(u^N)}{d u^N}$.
\item[(vi)] $F(u^N)$ is conserved under the flow of \eqref{PDE2}.
\end{enumerate}

\noi
Note that (v) holds true if the symplectic form $J$ commutes with the Dirichlet projection $\mathbb{P}_N$.
e.g. $J = i$ or $\dx$.
(vi) follows from (v) if $\mu$ is the Gibbs measure defined in \eqref{Gibbs}.
(vi) also follows easily if $F(u)$ consists only of the quadratic part such as $\int u^2$.
Just note that $\int u^N v^N = \int u^N v$.
Thus, in computing $\dt F(u^N)$ with \eqref{PDE2}, 
$\mathbb{P}_N$ in front of $\mathcal{N}(u^N)$
plays no role, and hence the computation reduces to that for \eqref{PDE1},
which is known to conserve $F$.

By Liouville's theorem and (v), 
the Lebesgue measure $\prod_{|n| \leq N} d \ft{u^N}(n)$ is invariant 
under the flow of \eqref{PDE2}.
Hence, 
the finite dimensional version $\mu_N$ of $\mu$:
\[d \mu_N = Z_N^{-1} e^{-F(u^N) }
\prod_{x \in \T} du^N(x)\]

\noi
is invariant under the flow of \eqref{PDE2}.
Assume that $\mu_N$ converges weakly to $\mu$. 
See Proposition \ref{PROP:Zhidkov2}.
Using the invariance of $\mu_N$, 
Bourgain \cite{BO4, BO6} proved the following estimate on $u^N$.

\begin{proposition} \label{PROP:BO}
Given $T > 0$ and  $\eps > 0$, there exists $ \Omega_N \subset B$
with $\mu_N (\Omega_N^c) < \eps$ 
such that 
for $u_0^N \in \Omega_N$, \eqref{PDE2} is well-posed on $[-T, T]$
with the growth estimate:
\[ \| u^N(t) \|_{B} \lesssim \bigg( \log \frac{T}{\eps} \bigg)^\frac{1}{2}, \text{ for }|t| \leq T.\]
\end{proposition}

\noi
In proving Proposition \ref{PROP:BO}, 
 we need to assume the following large deviation estimate.

\begin{lemma} \label{LEM:DEV}
There exists $c>0$, independent of $N$, such that for sufficiently large $K>0$, we have
\[ \mu_N \big( \{ \|u^N_0\|_B > K\}) < e^{-cK^2}.\]
\end{lemma}

\noi
Note that Lemma \ref{LEM:DEV} is basically 
Fernique's theorem \cite{FER}
since $(B, \mu)$ is an abstract Wiener space. 
See Theorem \ref{THM:FER} below.

\begin{proof}[Proof of Proposition \ref{PROP:BO}]
Let $\Phi_N(t) $ denote the flow map of  \eqref{PDE2},
and define
\[  \Omega_N = \bigcap_{j = -[T/\dl]}^{[T/\dl]} \Phi_N^j(\dl)(\{ \|u^N_0\|_B \leq K\}).\]

\noi
By invariance of $\mu_N$ and $\dl \sim K^{-\al}$, we have 
\[ \mu_N (\Omega_N^c) \lesssim \frac{T}{\dl} \mu_N (\{ \|u^N_0\|_B > K\})
\sim T K^\al e^{-cK^2}.\]

\noi
By choosing $K\sim \big( \log \frac{T}{\eps} \big)^\frac{1}{2}$, we have
$\mu_N (\Omega_N^c) < \eps$.
Moreover, 
by its construction, $\| u^N(j \dl) \|_B \leq K$ for $j = 0, \cdots, \pm [T/\dl]$.
By local theory, we have 
\[\| u^N(t) \|_B \leq 2K \sim \Big( \log \frac{T}{\eps} \Big)^\frac{1}{2} \text{ for } |t| \leq T.\]

\noi
Hence, $\Omega_N$ has the desired property.
\end{proof}

As a corollary to Proposition \ref{PROP:BO}, 
one needs to prove the following statements.

\begin{itemize}
\item[(a)]
Given $\eps > 0$, there exists $ \Omega_\eps \subset B$
with $\mu (\Omega_\eps^c) < \eps$ 
such that 
for $u_0 \in \Omega_\eps$, \eqref{PDE1} is globally well-posed 
with the growth estimate:
\begin{equation}\label{Gbound1}
 \| u(t) \|_{B} \lesssim \bigg( \log \frac{1+|t|}{\eps} \bigg)^\frac{1}{2}, \text{ for all } t \in \R.
\end{equation}

\item[(b)] \label{BB} The uniform convergence lemma:
\[ \|u - u^N\|_{C([-T, T];B')} \to 0\]
as $N \to \infty$ uniformly for $u_0\in\Omega_\eps$,
where $B' \supset B$.
\end{itemize} 

\noi  
Note that 
(a) implies that \eqref{PDE2} is a.s. GWP,
since $\wt{\Omega} := \bigcup_{\eps>0} \Omega_\eps$ has probability 1.
One can prove (a) and (b) 
by estimating the difference $u - u^N$ of 
solutions to 
\eqref{PDE1}  and  \eqref{PDE2},
using the estimates (ii)$\sim$(iv) and applying Proposition \ref{PROP:BO} to $u^N$.
We point out one useful observation due to Bourgain \cite{BO6}.
For KdV, the nonlinearity of the difference equation is given by
\[F(t) = \dx u^2(t) - \mathbb{P}_N \dx (u^N)^2(t).\]

\noi 
Since $\mathbb{P}_N \Big( \big( \mathbb{P}_\frac{N}{2}u\big)^2 \Big) 
=  \big( \mathbb{P}_\frac{N}{2}u\big)^2 $, we have
\begin{align} \label{DIFF}
F  =  \dx \Big(  u^2 - \big( \mathbb{P}_\frac{N}{2} u\big)^2 \Big)
+  \mathbb{P}_N \dx \Big( \big( \mathbb{P}_\frac{N}{2}u\big)^2 - u^2 \Big) 
+  \mathbb{P}_N\dx (u^2 - (u^N)^2).
\end{align}

\noi
After applying the nonlinear estimate, 
the first two terms can be made small due to the factor
$u -  \mathbb{P}_\frac{N}{2} u$,
and the last term has the factor $u - u^N$,
which we need to close the argument.

Finally, putting all the ingredients together, 
we obtain the invariance of $\mu$.
See the diagram below.
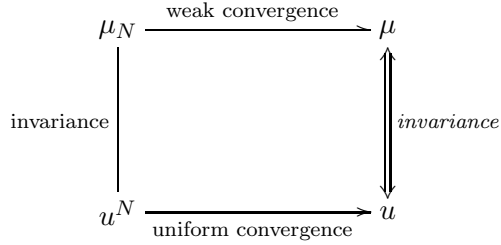
\begin{equation*} \xymatrix{\mu_N  \ar@{-}[dd]_{\text{invariance}} \ar@^{->}[rrr]^{\text{weak convergence}} 
& &  & \mu \ar@2{<->}[dd]^{\textit{invariance}}\\ \\ 
  u^N \ar[rrr]_{\text{uniform convergence}} & & & u  }\end{equation*}

 To conclude this subsection, 
we give several examples of the Banach spaces $B$ used for proving the invariance of  the Gibbs measures.
Note that for the radial results in $d = 2, 3$, 
$H^s$ denotes the Sobolev spaces in terms of the eigenfunctions of the Laplace operator on $\mathbb{D}^d$
(with appropriate boundary conditions.)

\smallskip
\noi
$\bullet$ $B = H^{\frac{1}{2}-}$:
 quintic or sub-quintic NLS, KdV \cite{BO4}, 
 subquintic radial NLS on $\mathbb{D}^2$ \cite{TZ1, TZ2},
 subquartic radial NLW on $\mathbb{D}^3$ with 
the Dirichlet boundary condition \cite{BT1, BT3}.

\smallskip

\noi
$\bullet$ $B = H^{\frac{1}{2}-} \cap \mathcal{F}L^{1-, \infty}$:
mKdV \cite{BO4}, Zakharov \cite{BO5}, coupled KdV systems with Diophantine conditions \cite{OH3}, 
Schr\"odinger-Benjamin-Ono \cite{OHSBO}.

\smallskip

As we saw already, we can not use $H^{-\frac{1}{2}-}$ 
to study the invariance of the white noise for KdV.
Hence, 
 we propose to use the Besov-type space $\ft{b}^s_{p, \infty}$, $sp <-1$,
 defined in  \eqref{Besov}.
In the following subsections, 
we show that $\ft{b}^s_{p, \infty}$ captures the low regularity of the white noise
for $sp < -1$,
and that KdV is LWP in  $\ft{b}^s_{p, \infty}$, $sp <-1$.

\subsection{Abstract Wiener spaces}

In Subsection 1.2, we reviewed the Gaussian measures in Hilbert spaces.
However, $\ft{b}^s_{p, \infty}$ is not a Hilbert space,
so  
we briefly go over  the basic theory of abstract Wiener spaces.

Recall the following definitions from Kuo \cite{KUO}:
Given  a real separable Hilbert space $H$ with norm $\|\cdot \|$, 
let $\mathcal{F} $ denote the set of finite dimensional orthogonal projections $\mathbb{P}$ of $H$.
Then, define a cylinder set $E$ by  $E = \{ x \in H: \mathbb{P}x \in F\}$ where $\mathbb{P} \in \mathcal{F}$ 
and $F$ is a Borel subset of $\mathbb{P}H$,
and let $\mathcal{R} $ denote the collection of such cylinder sets.
Note that $\mathcal{R}$ is a field but not a $\s$-field.
Then, the Gauss measure $\mu$ on $H$ is defined 
by 
\[ \mu(E) = (2\pi)^{-\frac{n}{2}} \int_F e^{-\frac{\|x\|^2}{2}} dx  \]

\noindent
for $E \in \mathcal{R}$, where
$n = \text{dim} \mathbb{P} H$ and  
$dx$ is the Lebesgue measure on $\mathbb{P}H$.
It is known that $\mu$ is finitely additive but not countably additive in $\mathcal{R}$.

A seminorm $|||\cdot|||$ in $H$ is called measurable if for every $\eps>0$, 
there exists $\mathbb{P}_\eps \in \mathcal{F}$ such that 
\[ \mu( ||| \mathbb{P} x ||| > \eps  )< \eps \]

\noindent
for $\mathbb{P} \in \mathcal{F}$ orthogonal to $\mathbb{P}_\eps$.
Any measurable seminorm  is weaker  than the norm of $H$,
and $H$ is not complete with respect to $|||\cdot|||$ unless $H$ is finite dimensional.
Let $B$ be the completion of $H$ with respect to $|||\cdot|||$
and denote by $i$ the inclusion map of $H$ into $B$.
The triple $(i, H, B)$ is called an abstract Wiener space.

Now, regarding $y \in B^\ast$ as an element of $H^\ast \equiv H$ by restriction,
we embed $B^\ast $ in $H$.
Define the extension of $\mu$ onto $B$ (which we still denote by $\mu$)
as follows.
For a Borel set $F \subset \R^n$, set
\[ \mu ( \{x \in B: ((x, y_1), \cdots, (x, y_n) )\in F\})
:= \mu ( \{x \in H: (\jb{x, y_1}_H, \cdots, \jb{x, y_n}_H )\in F\}),\]

\noindent
where $y_j$'s are in $B^\ast$ and $(\cdot , \cdot )$ denote the natural pairing between $B$ and $B^\ast$.
Let $\mathcal{R}_B$ denote the collection of cylinder sets
$ \{x \in B: ((x, y_1), \cdots, (x, y_n) )\in F \}$
in $B$.
Note that the pair $(B, \mu)$ is often referred to as an abstract Wiener space as well.

\begin{theorem}[Gross \cite{GROSS}]
$\mu $ is countably additive in the $\s$-field generated by $\mathcal{R}_B$.
\end{theorem}

\noi
In the present context, let $H= L^2(\mathbb{T})$ and
$B=\ft{b}^s_{p, \infty} (\mathbb{T})$ for $sp < -1$. 
Then, we have 

\begin{proposition} \label{PROP:meas}
The seminorm  $\|\cdot\|_{\ft{b}^s_{p, \infty}}$ is measurable for $sp < -1$.

\end{proposition}

\noindent
Hence, $(i, H, B) = (i, L^2, \ft{b}^s_{p, \infty}) $ is an abstract Wiener space, 
and $\mu$ defined in \eqref{White} is countably additive in $\ft{b}^s_{p, \infty}$.
We present the proof of Proposition \ref{PROP:meas} at the end of this subsection.
For our application, we can choose $s$ and $p$ such that $sp < -1$,
and thus we will not discuss the endpoint case. 
Also, note that in following Bourgain's argument as in Subsection 3.1, 
we need $sp < -1$ since we need a pair $(B, \mu)$, $(B', \mu)$ of abstract Wiener spaces
with $B \subset B'$. See (b) on p.\pageref{BB}.
It also follows from the proof that $(i, L^2, \mathcal{F} L^{s, p}) $, 
where $ \mathcal{F}L^{s, p} = \ft{b}^s_{p, p}$ defined in \eqref{FLP},
is also an abstract Wiener space for $sp < -1$ (we need  a strict inequality in this case.)

Given an abstract Wiener space $(i, H, B)$, we have the following integrability result due to Fernique \cite{FER}.

\begin{theorem} [Theorem 3.1 in \cite{KUO}]
\label{THM:FER} 
Let $(i, H, B)$ be an abstract Wiener space.
Then, there exists $ c > 0$ such that $ \int_B e^{c \|x\|_B^2} \mu(d x) < \infty$.
Hence, there exists $ c' > 0$ such that $\mu ( \|x\|_B > K) \leq e^{-c'K^2}$
for sufficiently large $K>0$.
\end{theorem}

\noi
In our context, if $sp < -1$, we have 
$ \mu \big(\|\phi\|_{\ft{b}^s_{p, \infty}(\mathbb{T})} \geq K, \phi \ \text{mean } 0 )  \leq e^{-c K^2}$ for some $c>0$.
With this estimate and Theorem \ref{THM:LWP2}, 
we can follow the argument in \cite{BO4} to prove Theorem \ref{THM:GWP2}.
We omit the details.
Also, see \cite{BT1}, \cite{OH3}, \cite{TZ1}, \cite{TZ2} for the details.

\begin{proof}[Proof of Proposition \ref{PROP:meas}]
We present the proof only for $2 < p < \infty$, which is the relevant case for our application.
We just point out that the proof for $p \leq 2$ is similar but simpler
(where one can use H\"older inequality in place of Lemma \ref{CL:decay} below.)
For $p = \infty$, see \cite{BO4}, \cite{BO5}, \cite{OH3}.

It suffices to show that
for given $\eps> 0$, 
there exists large $M_0$ such that 
\begin{equation*} 
 \mu \big(\|\mathbb{P}_{>M_0}\phi\|_{\ft{b}^s_{p, \infty}} > \eps) < \eps,
\end{equation*}

\noi 
where $\mathbb{P}_{>M_0}$is the projection onto the frequencies $|n| > M_0$.
In the following, 
write $\phi = \sum_{n\ne 0} g_n e^{inx}$, 
where $\{ g_n(\omega) \}_{n = 1}^\infty$ is a sequence of independent standard complex-valued Gaussian random variables
and $g_{-n} = \cj{g_n}$.
First, recall the following lemma.
\begin{lemma}[Lemma 4.7 in \cite{OHSBO}] \label{CL:decay}
Let $\{g_n\}$ be a sequence of i.i.d standard complex-valued Gaussian random variables.
Then, for $M$ dyadic and $\dl > 0$, we have
\[ \lim_{M\to \infty} M^{1-\dl} \frac{\max_{|n|\sim M } |g_n|^2}{ \sum_{|n|\sim M} |g_n|^2} = 0, \text{ a.s.} \]
\end{lemma}

Next, we present a large deviation lemma.
This can be proved by a direct computation using the polar coordinate.  See 
\cite{BO4}, \cite{OH3}, \cite{TZ1}.

\begin{lemma} \label{LEM:polar}

Let $M$ be dyadic, and $R = R(M)  \geq M^{\frac{1}{2}+} $.
Then, there exists $c$ such that
\begin{equation} \label{polar}
\mathbb{P}_\omega \big[\Big( \sum_{n\sim M} |g_n(\omega) |^2\Big)^\frac{1}{2} \geq R \big]
\leq  e^{-c R^2}
\end{equation}

\noi
for all dyadic $M$ (i.e. $c$ is independent of $M$.)
Moreover, this is essentially sharp in the sense that
\eqref{polar} can not hold if $R \leq M^\frac{1}{2}$.
\end{lemma}

Fix $K> 1$ and $\dl \in ( 0, \frac{1}{2})$ (to be chosen later.)
Then, by Lemma \ref{CL:decay} and Egoroff's theorem, 
there exists a set $E$ such that $\mu (E^c) < \frac{1}{2}\eps$
and the convergence in Lemma \ref{CL:decay} is uniform on $E$.
i.e. we can choose dyadic $M_0$ large enough such that 
\begin{equation} \label{A:decay}
\frac{\| \{g_n (\omega) \}_{|n| \sim M} \|_{L^{\infty}_n} }
{\| \{g_n (\omega) \}_{|n| \sim M} \|_{L^{2}_n} }
 \leq M^{-\dl}, 
\end{equation}

\noindent 
for  all $\omega \in E$  and dyadic $M > M_0$.
In the following, we will work only on $E$ and  drop `$\cap E$' for notational simplicity. 
However, it should be understood that all the events are under the intersection with $E$ 
so that  \eqref{A:decay} holds.

The basic idea of the following argument is due to Bourgain's dyadic pigeonhole principle
in \cite{BO4}.
Let $\{\s_j \}_{j \geq 1}$ be a sequence of positive numbers such that $\sum \s_j = 1$,
and let $ M_j = M_0 2^j$ dyadic.
Note that $\s_j = C 2^{-\ld j} =C M_0^\ld M_j^{-\ld}$ for some small $\ld > 0$
(to be determined later.)
Then,  we have
\begin{align} \label{A:subadd} 
\mu \big(\|\mathbb{P}_{>M_0}\phi\|_{\ft{b}^s_{p, \infty}} > \eps)
& \leq \mu  \big( \| \{g_n \}_{|n| > M_0} \|_{{b}^s_{p, 1}} > \eps \big) \notag \\
& \leq \sum_{j = 0}^\infty \mu \big( \| \{\jb{n}^s g_n \}_{|n| \sim M_j} \|_{L_n^{p}}  > \s_j \eps \big),
\end{align}

\noindent
where ${b}^s_{p, 1}$ is as in \eqref{Besov}
with the $l^\infty$ norm over the dyadic blocks replaced by 
the $l^1$ sum.
By interpolation and \eqref{A:decay}, we have  
\begin{align*}
\| \{ & \jb{n}^s g_n \}_{|n| \sim M_j} \|_{L_n^{p}} 
\sim M_j^{s} \| \{ g_n \}_{|n| \sim M_j} \|_{L_n^{p}} 
\leq  M_j^{s} \| \{ g_n \}_{|n| \sim M_j} \|_{L_n^{2}}^\frac{2}{p}
  \| \{ g_n \}_{|n| \sim M_j} \|_{L_n^{\infty}}^\frac{p-2}{p} \\
& \leq  M_j^{s} \| \{ g_n \}_{|n| \sim M} \|_{L_n^{2}}
\Bigg(\frac{  \| \{ g_n \}_{|n| \sim M_j} \|_{L_n^{\infty}} }
{\| \{ g_n \}_{|n| \sim M_j} \|_{L_n^{2}}} \Bigg)^\frac{p-2}{p}
\leq  M_j^{s -\dl \frac{p-2}{p}} \| \{ g_n \}_{|n| \sim M_j} \|_{L_n^{2}}
\end{align*}

\noindent
a. s. 
Thus, if we have $\|\{\jb{n}^s g_n \}_{|n| \sim M_j} \|_{L_n^{p}}  > \s_j \eps$,
 then we have
 $\| \{ g_n \}_{|n| \sim M_j} \|_{L_n^{2}} 
 \gtrsim R_j $
 where $R_j := \s_j \eps M_j^{-s+\dl \frac{p-2}{p}} $.
With $p = 2 + 2\theta$, we have 
$-s+\dl \frac{p-2}{p} = \frac{-sp + 2 \dl \theta}{2 + 2 \theta} > \frac{1}{2}$
by taking $\dl$ sufficiently close to $\frac{1}{2}$ since $-sp > 1$.
Then, by taking $\ld > 0$ sufficiently small,
$R_j = \s_j \eps M_j^{-s+\dl \frac{p-2}{p}} 
= C \eps M_0^\ld M_j^{-s+\dl \frac{p-2}{p}-\ld}
\gtrsim C \eps M_0^{\ld} M_j^{\frac{1}{2}+} $. 
Then, by Lemma \ref{LEM:polar},  we have
\begin{align} \label{A:highfreq1}
\mu\big(   \| \{ g_n \}_{|n| \sim M_j} \|_{L_n^{2}}  \gtrsim R_j \big)
\leq C \int_{R_j}^\infty e^{-\frac{1}{4}r^2} r dr  
\leq e^{-cR_j^2} = e^{-cC^2 M_0^{2\ld} M_j^{1+} \eps^2}.
\end{align}

\noindent
From \eqref{A:subadd} and \eqref{A:highfreq1}, we have
\begin{align*} 
\mu \big(\|\mathbb{P}_{>M_0}\phi\|_{\ft{b}^s_{p, \infty}} > \eps)
\leq \sum_{j =1}^\infty 
e^{-cC^2 M_0^{1+2\ld+} 2^{j+} \eps^2} \leq \tfrac{1}{2} \eps
\end{align*}

\noi
by choosing $M_0$ sufficiently large.
\end{proof}

\subsection{Function spaces and basic  embeddings}

First, let $X^{s, b}$ denote the usual periodic Bourgain space defined in \eqref{Xsb}.
We often use the shorthand notation $\|\cdot\|_{s, b}$ to denote the $X^{s, b}$ norm.
Now, define $X^{s, b}_{p, q}$, the Bourgain space adapted to $\ft{b}^s_{p, \infty}$, 
 to be the completion of the Schwartz functions on $\T\times \R$
with respect to the norm given by
\begin{equation} \label{XSBPQ}
 \| u \|_{X^{s, b}_{p, q}} 
 =  \|\jb{n}^s \jb{\tau - n^3}^b \ft{u}(n, \tau)\|_{b^0_{p, \infty} L^q_\tau}
 = \sup_j \|\jb{n}^s  \jb{\tau - n^3}^b \ft{u}(n, \tau)
\|_{L^p_{|n| \sim 2^j} L^q_\tau}.
\end{equation}

\noi
In the following, we take $p = 2+$ and $s = -\frac{1}{2}+ = -\frac{1}{2}+\dl$  
with $\dl < \frac{p-2}{2p}$ (and $\dl > \frac{p-2}{4p}$)
such that $sp < -1$.
Lastly, given $T> 0$, we define $X^{s, b, T}_{p, q}$ as a restriction of $X^{s, b}_{p, q}$ on $[0, T]$
by
\[ \|u\|_{X^{s, b, T}_{p, q}} 
=  \|u\|_{X^{s, b}_{p, q}[0, T]} 
= \inf \big\{ \|\wt{u} \|_{X^{s, b}_{p, q}}: {\wt{u}|_{[0, T]} = u}\big\}.\]

\noi
We define the local-in-time versions of other function spaces analogously.

Now, we discuss the basic embeddings.
For $p \geq 2$, we have $\| a_n\|_{L^p_n} \leq \| a_n\|_{L^2_n}$.
Thus, we have $\|f\|_{\ft{b}^s_{p, \infty}} \leq \|f\|_{H^{s}}$, and thus
\begin{equation}\label{EMBED1}
\|u\|_{X^{s, b}_{p, 2}} \leq \|u\|_{X^{s, b}}.
\end{equation}

\noi
By H\"older inequality, we have
\begin{align} \label{EMBED2}
\|f\|_{H^{-\frac{1}{2}-\dl}} 
& \leq \sup_j \| \jb{n}^{-2\dl+} \|_{L^\frac{2p}{p-2}} \|\jb{n}^{-\frac{1}{2} + \dl} \ft{f}(n) \|_{L^p_n}
\leq \| f \|_{\ft{b}^s_{p, \infty}} 
\end{align}

\noi
for $s = -\frac{1}{2} +\dl$ with $\dl > \frac{p-2}{4p}$.
Hence, for $s = - \frac{1}{2} +\dl$ with $\dl > \frac{p-2}{4p}$, 
we have
\begin{equation} \label{EMBED3}
\|u\|_{X^{-\frac{1}{2}-\dl, b}}\lesssim \|u\|_{X_{p, 2}^{s, b}}.
\end{equation}

Now, we briefly go over the linear  estimates.
Let $S(t) = e^{-t \dx^3}$ and $T\leq 1$ in the following.
We first present the homogeneous and nonhomogeneous linear estimates.
See \cite{BO1}, \cite{OH4} for details.

\begin{lemma} \label{LEM:linear1}
For any $s \in \mathbb{R}$ and $b < \frac{1}{2}$, we have 
$\|  S(t) u_0\|_{X^{s, b, T}_{p, 2}} \lesssim T^{\frac{1}{2}-b}\|u_0\|_{\ft{b}^s_{p, \infty}}$.
\end{lemma}

\begin{lemma} \label{LEM:linear2}
For any $s \in \mathbb{R}$ and $b \leq \frac{1}{2}$, we have 
\begin{align*} 
 \bigg\|  \int_0^t S(t-t') F(x, t') dt'\bigg\|_{X^{s, b, T}_{p, 2}} 
\lesssim \| F \|_{X^{s, b-1}_{p, 2}} + \| F \|_{X^{s, -1}_{p, 1}}.
\end{align*}

\noi
Also, we have 
$ \Big\|  \int_0^t S(t-t') F(x, t') dt'\Big\|_{X^{s, b, T}_{p, 2}} 
\lesssim \| F \|_{X^{s, b-1}_{p, 2}}$
for $b > \frac{1}{2}$.
\end{lemma}

\noi
The next lemma is the periodic $L^4$ Strichartz estimate due to Bourgain \cite{BO1}.
\begin{lemma} \label{LEM:L4}
Let $u$ be a function on $\T \times \R$.
Then, we have 
$ \|u\|_{L^4_{x, t}} \lesssim \|u\|_{X^{0, \frac{1}{3}}}.$
\end{lemma}

\subsection{Nonlinear analysis}

Now, we present the crucial nonlinear analysis.
First, we briefly go over Bourgain's argument in \cite{BO3}.
By writing the integral equation,  KdV \eqref{KDV} is equivalent to
\begin{equation} \label{KDVduhamel}
u(t) =  S(t) u_0 -\tfrac{1}{2} \mathcal{N}(u, u) (t), 
\end{equation} 

\noi
where $\mathcal{N}(\cdot, \cdot)$  is  given by
\begin{equation} \label{NN}
\mathcal{N}(u_1, u_2) (t) :=  \int_0^t S(t - t') \dx (u_1 u_2)(t') d t'.
\end{equation}

In the following, we assume that the initial condition $u_0$ has the mean 0,
which implies that $u(t)$ has the spatial mean 0 for each $t\in \R$.
We use $(n, \tau)$, $(n_1, \tau_1)$, and $(n_2, \tau_2)$ to denote the Fourier variables
for $uu$, the first factor, and the second factor $u$ of $uu$ in $\mathcal{N}(u, u)$, 
respectively.
i.e. we have $n = n_1 + n_2$ and $\tau = \tau_1 + \tau_2$.
By the mean 0 assumption on $u$ and 
by the fact that we have $\dx (uu)$ in the definition of $\mathcal{N}(u,u)$, 
we assume  $n, n_1, n_2 \ne 0$.
We also use the following notation:
\[\s_0 := \jb{\tau - n^3} \text{ and } \s_j := \jb{\tau_j - n_j^3}.\]

\noi
One of the main ingredients is the observation due to Bourgain \cite{BO1}:
\begin{equation} \label{Walgebra}
n^3 - n_1^3 - n_2^3 = 3 n n_1 n_2, \ \text{for } n = n_1 + n_2,
\end{equation}

\noi
which in turn implies that 
\begin{equation} \label{MAXMAX}
\MAX:= \max( \s_0, \s_1, \s_2) \gtrsim \jb{n n_1 n_2}.
\end{equation}

Now, define 
\begin{equation} \label{AJJ}
A_j = \{(n, n_1, n_2, \tau, \tau_1, \tau_2) \in \mathbb{Z}^3 \times \R^3:
\s_j = \MAX\},
\end{equation}

\noi
and let $\mathcal{N}_j(u, u)$ denote the contribution of $\mathcal{N}(u, u)$ on $A_j$.
By the standard bilinear estimate with Lemma \ref{LEM:L4} as in \cite{BO1}, \cite{KPV4}, 
we have
\begin{align} \label{N_0}
\|\mathcal{N}_0(u, u)\|_{{-\frac{1}{2} + \dl, \frac{1}{2}-\dl}}
\leq o(1)\|u\|^2_{{-\frac{1}{2} - \dl, \frac{1}{2}-\dl}},
\end{align}

\noi
where $o(1) = T^\theta$ with $\theta > 0$ by considering the estimate 
on a short time interval $[-T, T]$.
See (2.17), (2.26), and (2.68) in \cite{BO3}.
Here, we abuse the notation and use $\|\cdot\|_{s, b} = \|\cdot\|_{X^{s, b}}$ 
to denote the local-in-time version as well.
Note that the temporal regularity is $b = \frac{1}{2} - \dl < \frac{1}{2}$.
This allowed us to gain the spatial regularity by $2\dl$.
Clearly, we can not expect to do the same for $\mathcal{N}_1(u, u)$.
(By symmetry, we do not consider $\mathcal{N}_2(u, u)$ in the following.)
The bilinear estimate \eqref{KPVbilinear} is known to fail for any $s \in \R$
if $b < \frac{1}{2}$ due to the contribution from $\mathcal{N}_1(u, u)$. 
See \cite{KPV4}.
Following the notation in \cite{BO3}, 
let 
\begin{equation}\label{Isb}
I_{s, b} = \|\mathcal{N}_1(u, u) \|_{X^{s, b}} \
\text{ and }\ \al := \frac{1}{2} -\dl < \frac{1}{2}.
\end{equation}

\noi
Then, by Lemma \ref{LEM:linear2} and duality 
with $\|d(n, \tau)\|_{L^2_{n, \tau}} \leq 1$, we have 
\begin{align} \label{eq:I1}
I_{-\al, 1-\al} & = \|\mathcal{N}_1(u, u) \|_{-\al, 1-\al}\\
& \lesssim \sum_{\substack{n, n_1\\n = n_1 + n_2}} \intt_{\tau = \tau_1 + \tau_2} d\tau d\tau_1
\frac{\jb{n}^{1-\al}d(n, \tau)}{\s_0^{\al }} \ft{u}(n_1, \tau_1) \frac{\jb{n_2}^{1-\al}c(n_2, \tau_2)}{\s_2^\al},
\notag
\end{align}

\noi
where  
\begin{equation} \label{CN}
c(n_2, \tau_2) = \jb{n_2}^{-(1-\al)}\s_2^{\al}\, \ft{u}(n_2, \tau_2) 
\text{ so that }
\|c\|_{L^2_{n, \tau}} = \|u\|_{-(1-\al), \al} = \|u\|_{-\frac{1}{2}-\dl, \frac{1}{2}-\dl}.
\end{equation}

\noi
The main idea here is to consider the second iteration, 
i.e. substitute \eqref{KDVduhamel} for  $\ft{u}(n_1, \tau_1)$ in \eqref{eq:I1}, 
thus leading to a trilinear expression.
Since $\s_1 = \MAX \gtrsim \jb{n n_1n_2}\gg1$ on $A_1$, 
we can assume that 
\begin{equation} \label{eq:u_1}
\ft{u}(n_1, \tau_1) = \big(\mathcal{N}(u, u)\big)^\wedge(n_1, \tau_1)
\sim \frac{|n_1|}{\s_1} \sum_{n_1 = n_3 + n_4} \intt_{\tau_1 = \tau_3 + \tau_4}
\ft{u}(n_3, \tau_3)\ft{u}(n_4, \tau_4) d\tau_4.
\end{equation}

\noi
Note that the $\s_1 $ appearing in the denominator
allows us to cancel $\jb{n}^{1-\al}$ and $\jb{n_2}^{1-\al}$ in the numerator in \eqref{eq:I1}.
Then, $I_{-\al, 1-\al}$ can be estimated by 
\begin{align} \label{eq:I2}
\lesssim \sum_{\substack{n = n_1 + n_2\\n_1 = n_3 + n_4}} 
\intt_{\substack{\tau = \tau_1 + \tau_2\\\tau_1 = \tau_3 + \tau_4}} 
\frac{\jb{n}^{1-\al}d(n, \tau)}{\s_0^{\al }} \frac{|n_1|}{\s_1} \, \ft{u}(n_3, \tau_3)\ft{u}(n_4, \tau_4)
\frac{\jb{n_2}^{1-\al}c(n_2, \tau_2)}{\s_2^\al}
.
\end{align}

\noi
The argument was then divided into several cases, 
depending on the sizes of $\s_0, \cdots, \s_4$.
Here, the key algebraic relation is
\begin{equation}\label{algebra2}
n^3 - n_2^3 - n_3^3 - n_4^3 = 3(n_2+ n_3)(n_3+ n_4)(n_4+ n_2), \  \text{ with } n = n_2 + n_3 + n_4.
\end{equation}

\noi
Then, Bourgain proved -see (2.69) in \cite{BO3}-
\begin{equation} \label{N_1}
I_{-\al, 1-\al} \leq o(1)\|u\|_{-(1-\al), \al} I_{-\al, 1-\al}  
+ o(1) \|u\|^3_{-(1-\al), \al} + o(1) \|u\|_{-(1-\al), \al},
\end{equation}

\noi
{\it assuming} the a priori estimate \eqref{BOO}: $|\ft{u}(n, t)| <C$ for all $n\in \mathbb{Z}$, $t \in\R$.
Indeed, the estimates involving the first two terms on the right hand side of  \eqref{N_1}
were obtained without \eqref{BOO}, and
{\it only} the last term in \eqref{N_1} required \eqref{BOO}, 
-see ``Estimation of (2.62)'' in \cite{BO3}-, 
which was then used to deduce
\begin{equation} \label{eq:apriori}
\|\ft{u}(n, \cdot)\|_{L^2_\tau} < C.
\end{equation}

\noi
The a priori estimate \eqref{BOO} is derived via
the isospectral property of the KdV flow
and is false for a general function in $X^{-(1-\al), \al}$.
(It is here that the smallness of the total variation $\|\mu\|$ is used.)

\medskip

Our goal is to carry out a similar analysis on the second iteration 
{\it without} the a priori estimates \eqref{BOO} and \eqref{eq:apriori} coming from the complete integrability of KdV.
We achieve this goal by considering the estimate
in $X^{-\al, \al}_{p, 2} = X^{-\frac{1}{2}+ \dl, \frac{1}{2}-\dl}_{p, 2}$, 
where $p = 2+$ and $\frac{p-2}{4p} < \dl< \frac{p-2}{2p}$.
By \eqref{EMBED1} and \eqref{EMBED3}
(recall $-\al = -\frac{1}{2} +\dl$ and $-(1-\al) = -\frac{1}{2} -\dl$), we have
\begin{equation} \label{EMBED5}
\|u\|_{X^{-\al, \al}_{p, 2}} \leq \|u\|_{X^{-\al, \al}},
\text{ and } \
\|u\|_{X^{-(1-\al), \al}} \lesssim \|u\|_{X^{-\al, \al}_{p, 2}}.
\end{equation}

\noi
Then, it follows from \eqref{N_0} and \eqref{EMBED5} that 
\begin{equation} \label{Np0}
\|\mathcal{N}_0(u, u)\|_{X^{-\al, \al}_{p, 2}} \leq o(1) \|u\|^2_{X^{-\al, \al}_{p, 2}}.
\end{equation}

Now, we consider the estimate on $\|\mathcal{N}_1(u, u) \|_{X^{-\al, \al}_{p, 2}}$.
From \eqref{EMBED5} and $\al < 1-\al$, it suffices to control $I_{-\al, 1-\al}$.
As before, we consider the second iteration, 
and substitute \eqref{KDVduhamel} for  $\ft{u}(n_1, \tau_1)$ in \eqref{eq:I1}.
Note that we can use the estimates on 
 $\mathcal{N}_1(\mathcal{N}(u, u), u)$ from \cite{BO3} 
{\it except} when the a priori bound \eqref{BOO} was assumed.
i.e. we need to estimate the contribution from (2.62) in \cite{BO3}:
\begin{equation}\label{262}
R_\al := 
\sum_{n} \intt_{\tau = \tau_2 +  \tau_3 + \tau_4} 
\chi_B \frac{d(n, \tau)}{\jb{n}^{1+\al}\s_0^{\al }} 
 \ft{u}(-n, \tau_2) \ft{u}(n, \tau_3)\ft{u}(n, \tau_4)d\tau_2 d\tau_3d\tau_4,
\end{equation}

\noi
where $\|d(n, \tau)\|_{L^2_{n, \tau}} \leq 1$
and $B = \{ \s_0, \s_2, \s_3, \s_4 < |n|^\g\}$ with some small parameter $\g>0$. 
Note that this corresponds to the case $n_2 = -n$ and $n_3 = n_4 = n$
in \eqref{eq:I2} after some reduction.
In our analysis, we directly estimate $R_\al$ 
in terms of $\|u\|_{X^{-\al, \al}_{p, 2}}$.
The key observation is that 
we can take the spatial regularity $s = -\al $ to be greater than $-\frac{1}{2}$ 
by choosing $p > 2$.

In \cite{BO3}, the parameter $\g = \g(\al)$,  subject to the conditions (2.43) and (2.60) in \cite{BO3},  
played a certain role in estimating $R_\al$
along with the a priori bound \eqref{BOO}.
However, it plays no role in our analysis.
By Cauchy-Schwarz and Young's inequalities, we have
\begin{align}
\eqref{262} & \leq \sum_n \| d(n, \cdot)\|_{L^2_\tau}
\jb{n}^{-1-\al} \| \ft{u}(-n, \tau_2)\|_{L^{\frac{6}{5}}_{\tau_2}}
\| \ft{u}(n, \tau_3)\|_{L^{\frac{6}{5}}_{\tau_3}}
\| \ft{u}(n, \tau_4)\|_{L^{\frac{6}{5}}_{\tau_4}} \notag
\intertext{By H\"older inequality (with appropriate $\pm$ signs) and the fact that $-1-\al < -3\al$,}
 & \leq \sum_n \| d(n, \cdot)\|_{L^2_\tau}
\prod_{j = 2}^4 \jb{n}^{-\al-}\|\s_j^{-\al}\|_{L^3_{\tau_j}}
\| \s_j^{\al} \ft{u}(\pm n, \tau_j)\|_{L^{2}_{\tau_j}} \label{RA}\\
 & \leq \| d(\cdot, \cdot)\|_{L^2_{n, \tau}}
\|u\|_{X^{-\al, \al}_{6, 2}}^3
\leq \|u\|_{X^{-\al, \al}_{p, 2}}^3, \notag
\end{align}

\noi
where the last two inequalities follow by choosing  $\al > \frac{1}{3}$ and $p = 2+ <6$.

\medskip

Now, we put all the {\it a priori} estimates together.
Fix $\al = \frac{1}{2} - \dl$ as in \eqref{Isb}.
From Lemma \ref{LEM:linear1}, we have 
\begin{equation} \label{WT0}
\| S(t) u_0\|_{X^{s, b, T}_{p, 2}} \leq C_1 \|u_0\|_{\ft{b}^{s}_{p, \infty}}
\end{equation}

\noi
for any $s, b\in \R$ with $C_1 =  C_1(b)$.
From the definition of $\mathcal{N}_j(\cdot, \cdot)$ and \eqref{Np0},
we have
\begin{align}  \label{WT1}
\|\mathcal{N}(u, u)\|_{X^{-\al, \al, T}_{p, 2}}
\leq C_2 T^{\theta_1} \|u\|^2_{X^{-\al, \al, T}_{p, 2}}
+ 2 \|\mathcal{N}_1(u, u)\|_{X^{-\al, \al, T}_{p, 2}}.
\end{align}

\noi
From \eqref{Isb} and \eqref{EMBED5}, we have
\begin{align} \label{WT2}
\|\mathcal{N}_1(u, u)\|_{X^{-\al, 1- \al, T}_{p, 2}} \leq I_{-\al, 1-\al}.
\end{align}

\noi
Also, from \eqref{N_1} and \eqref{RA}, we have 
\begin{align*}
I_{-\al, 1-\al} \leq C_3\big( T^{\theta_2} \|u\|_{X^{-\al, \al, T}_{p, 2}}  I_{-\al, 1-\al}  
+ T^{\theta_3} \|u\|_{X^{-\al, \al, T}_{p, 2}}^3 
\big).
 \end{align*}
 
For fixed $R > 0$, 
choose $T>0$ small such that $C_3T^{\theta_2} R \leq \frac{1}{2}$.
Then, we have 
\begin{align} \label{WT4}
I_{-\al, 1-\al} \leq 2C_3 
T^{\theta_3} \|u\|_{X^{-\al, \al, T}_{p, 2}}^3,
 \end{align}

\noi
for $\|u^N\|_{X^{-\al, \al, T}_{p, 2}} \leq R$.

Let $u$ and $v$ be the solutions of \eqref{KDV} with the initial data $u_0$ and $v_0$.
Then, 
from \eqref{WT0}$\sim$\eqref{WT4}, we have
\begin{align} \label{GAMMA1}
\| u \|_{X^{-\al, \al, T}_{p, 2}}
   \leq C_1 \|u_0\|_{\ft{b}^{-\al}_{p, \infty}}
+ \tfrac{1}{2} C_2 T^{\theta_1} \|u\|^2_{X^{-\al, \al, T}_{p, 2}} + 2C_3 
T^{\theta_3} \|u\|_{X^{-\al, \al, T}_{p, 2}}^3 
\end{align}

\noi
and
\begin{align} \label{GAMMA2}
\| u -  v\|_{X^{-\al, \al, T}_{p, 2}}  
  & \leq  C_1 \|u_0- v_0\|_{\ft{b}^{-\al}_{p, \infty}}
+  \tfrac{1}{2} C_2 T^{\theta_1} (\|u\|_{X^{-\al, \al, T}_{p, 2}} 
+\|v\|_{X^{-\al, \al, T}_{p, 2}})\|u - v\|_{X^{-\al, \al, T}_{p, 2}}\notag  \\
& \hphantom{XX}
+ C_5 T^{\theta_3} 
\big(\|u\|_{X^{-\al, \al, T}_{p, 2}}^2 + \|v\|_{X^{-\al, \al, T}_{p, 2}}^2\big) 
 \|u -v \|_{X^{-\al, \al, T}_{p, 2}} .
\end{align}

Note that in estimating the difference $u-v$ on $A_1$,  one needs to consider 
\begin{equation}\label{GAMMA3}
 \wt{I}_{-\al, 1-\al} := \| \mathcal{N}_1(u, u) - \mathcal{N}_1(v, v)\|_{-\al, 1-\al} 
\end{equation}

\noi
as in \cite{BO3}.
We can follow the argument on pp.135-136 in \cite{BO3},
except for $R_\al$ defined in \eqref{262}, yielding the third term in \eqref{GAMMA2}. 
As for $R_\al$, we can write 
\begin{align} \label{WTN1}
\mathcal{N}(\mathcal{N}(u,u), u)- \mathcal{N}(\mathcal{N}(v,v), v) 
= \mathcal{N}(\mathcal{N}(u+v,u-v),u)
+ \mathcal{N}(\mathcal{N}(v,v), u - v)
\end{align}

\noi as in (3.4) in \cite{BO3},
and then we can repeat the computation done for $R_\al$, 
also yielding the third  term in \eqref{GAMMA2}.

Now, we are basically done.
Fix mean zero  $u_0 \in \ft{b}^{-\al}_{p, \infty}(\T)$, 
and take a sequence of smooth $u_0^{(n)}$ converging to $u_0$ in 
$\ft{b}^{-\al}_{p, \infty}(\T)$.
Strictly speaking, one needs to assume that $u_0$ is in a slightly 
more regular space due to the $L^\infty$-nature of the norm. See \cite{OH6}.
Choose $T$ sufficiently small such that
$C_3 T^{\theta_2}R \leq \frac{1}{2}$,
where 
$R= 2C_1 \|u_0\|_{\ft{b}^{-\al}_{p, \infty}} + 1$.
Then, we see that the smooth global solutions $u^{(n)}$ with the initial data  $u_0^{(n)}$
converge in $X^{-\al, \al, T}_{p, 2}$.
Let $u$ denote the limit.
We still need to show
\begin{itemize}
\item[(i)] $u$ is indeed a solution to \eqref{KDV} with $u_0 \in \ft{b}^{-\al}_{p, \infty}(\T)$
as its initial datum.
\item[(ii)] $u \in C([-T, T]; \ft{b}^{-\al}_{p, \infty})$.
\item[(iii)] the uniqueness of solution and the continuous dependence on initial data.
\end{itemize}

\noi
The argument here is just a little extension of what we have done, so
we omit the details. See \cite{OH6}.
We point out that 
 the solution map is H\"older continuous with the bound
\[\|u-v  \|_{C([-T, T]; \ft{b}^{-\al}_{p, \infty})}
\leq C_1(R)\|u_0-v_0\|^\beta_{\ft{b}^{-\al}_{p, \infty}}\]

\noi
for some small $\beta > 0$. 
Nonetheless, we have
\begin{equation} \label{GOODESTI}
\|u-v  \|_{C([-T, T]; \ft{b}^{-\al}_{p, \infty})}
\leq C_2(R)\|u_0-v_0\|_{\ft{b}^{-\al}_{p, \infty}} 
+ C_3(R) \|u - v\|_{X^{-\al, \al, T}_{p, 2}}.
\end{equation}

\noi
Note that 
\eqref{GOODESTI} is a ``good'' estimate which lets us proceed 
with the approximation argument described in Subsection 3.1
to establish a.s. GWP and the invariance of the white noise.

\subsection{Application: Stochastic KdV}

As an application of the nonlinear estimate in Subsection 3.4, 
we present sharp local well-posedness of the periodic stochastic KdV equation (SKdV)
with additive space-time white noise:
\begin{equation} \label{SKDV}
\begin{cases}
du + (\dx^3 u + u\dx u) dt = \phi dW\\
u(x, 0) = u_0(x), 
\end{cases}
\end{equation}

\noindent
where $\phi$ is a bounded linear operator on $L^2(\T)$
and 
$W (t) = \frac{\partial B}{\partial x}$ is a cylindrical Wiener process on $L^2(\T)$.
i.e. 
$W (t) =  \beta_0(t) e_0 + \sum_{n \ne 0 } \frac{1}{\sqrt{2}} \beta_n(t) e_n(x)$
with $e_n(x) = \frac{1}{\sqrt{2\pi}}e^{inx}$, 
where $\{ \beta_n\}_{n \geq 0}$ is a family of mutually independent complex-valued Brownian motions
(here we take $\beta_0$ to be real-valued)
in a fixed probability space $(\Omega, \mathcal{F}, P)$
associated with a filtration $\{\mathcal{F}_t\}_{t \geq 0}$
and $\beta_{-n}(t) = \cj{\beta_n}(t)$ for $n \geq 1$.

In \cite{DDT1}, 
de Bouard-Debussche-Tsutsumi established LWP in $B^\s_{2, 1}$ with $s > \s \geq -\frac{1}{2}$ 
and $\phi$ is Hilbert-Schmidt from $L^2$ to $H^{{s}}$.
Their argument was based on the result by Roynette \cite{ROY}
on the endpoint regularity of the Brownian motion.
i.e. the Brownian motion $\beta(t)$ belongs to the Besov space $B^{1/2}_{p, q}$
if and only if $ q= \infty$ (with $1\leq p < \infty$.)
Then, they proved a variant of the bilinear estimate \eqref{KPVbilinear}
by Kenig-Ponce-Vega
adjusted to their Besov space setting,
establishing LWP via the fixed point theorem.
Note that the use of a variant of \eqref{KPVbilinear}
required a slight regularization of the noise in space via $\phi$
so that the smoothed noise has the spatial regularity $s >-\frac{1}{2}$.
Thus, they could not treat the space-time white noise, i.e. $\phi =$ Id,
which is Hilbert-Schmidt from $L^2$ to $H^{{s}}$ for $s < -\frac{1}{2}$.

Now, observe that $\beta(t)$ has a Gaussian distribution for each fixed $t$.
Thus, $W(t)$ has the same regularity as the spatial white noise for each $t$.
Also, recall that 
$\mathcal{F} L^{b-1, p}_\tau$ captures the (temporal) regularity of the white noise
if 
$(b-1) \cdot p < -1.$
Hence, 
our Bourgain space $X^{s, b}_{p, q}$ in \eqref{XSBPQ} captures the regularity of the  space-time white noise
for $sp < -1$, $b  <  \frac{1}{2}$, and $q = 2$.
We can indeed control the {\it stochastic convolution}:
\begin{equation} \label {stoconv}
\Phi(t) = \int_0^t S(t - t') dW(t') 
\end{equation}

\noi
appearing in the Duhamel formulation of \eqref{SKDV} with the estimate:
\[ \mathbb{E} \big( \|  \Phi \|_{X^{s, \frac{1}{2}-}_{p, 2} [0, T]} \big) \lesssim C(\eta, s, p) < \infty\]

\noi
for $s p<-1$.
See Proposition 4.1 in \cite{OH6}.
Combining this with the nonlinear estimate on the second iteration, we obtain:
\begin{theorem}
The stochastic KdV \eqref{SKDV} with additive space-time white noise, 
i.e. $\phi = \textup{Id}$,
is locally well-posed almost surely
(with the prescribed mean on $u_0$.)
\end{theorem}

Several remarks are in order.
In the nonlinear analysis on the second iteration, 
we have an extra term arising from the stochastic convolution \eqref{stoconv}.
i.e. we need to estimate $\mathcal{N}(\Phi, u)$. 
Thus, our nonlinear analysis is stochastic,
whereas the bilinear estimate in \cite{DDT1} is entirely deterministic.
Moreover, there is {\it no}
smooth solutions for \eqref{SKDV} with $\phi = \text{Id}$.
Hence, we need to construct smooth approximating solutions $u^{(n)}$
with smooth initial data $u^{(n)}_0$
as well as smooth covariance operators $\phi^{(n)}$.
See \cite{OH6} for details.

\section{Method 2: Probabilistic approach}

\subsection{General framework}

First, we briefly discuss the basic structure of the argument
for the $p = 4$ case.
In order to prove Theorem \ref{THM:WC2},
it suffices to show that, for any smooth mean 0 function $f$ on $\T$,  
\begin{equation} \label{mainconv}
C_\beta  \int  e^{i \int f u + \beta \int u^4 } \chi_{\{ \int u^2 \leq K \beta^{-\frac{1}{2}}\}} d \mu_\beta 
\to  e^{-\frac{1}{2} \|f\|_{L^2}^2},
\end{equation}

\noindent
for some $C_\beta$
where $\mu_\beta$ is defined in \eqref{Gauss2}, since \eqref{mainconv} implies
\[ \int e^{i \int f u} d \rho_\beta^{(4)} =  
\frac{C_\beta \int  e^{i \int f u + \beta \int u^4 } 
\chi_{\{ \int u^2 \leq K \beta^{-\frac{1}{2}}\}} d \mu_\beta }
{C_\beta \int  e^{ \beta \int u^4 } 
\chi_{\{ \int u^2 \leq K \beta^{-\frac{1}{2}}\}} d \mu_\beta } 
\to  \frac{ e^{-\frac{1}{2} \|f\|_{L^2}^2}}{e^{-\frac{1}{2} \|0\|_{L^2}^2}}
= e^{-\frac{1}{2} \|f\|_{L^2}^2}.\]

In order to show \eqref{mainconv}, 
we divide the space into several regions
depending on the variations of $\int u^2$ and $\int u^4$.
For this purpose, we introduce the $\beta$-Wick ordered polynomials:
\begin{align} \label{Wick2}
 : u^2 \! :_\beta & = u^2 -  a_\beta,\\
 \label{Wick4}
 : u^4 \! :_\beta & = u^4 - 6 a_\beta u^2 + 3 a_\beta^2,
\end{align}

\noindent
where
$ a_\beta = \mathbb{E}_{\mu_\beta}\big[\int u^2\big]= \sum_{n \ne 0 } \frac{1}{1 + \beta n^2}, $
and $\mathbb{E}_{\mu_\beta}$ denotes the expectation with respect to 
$\mu_\beta$.
Note that $a_\beta \sim  \pi \beta^{-\frac{1}{2}}$
from a Riemann sum approximation.
Then, by direct computation, we can show the following.

\begin{lemma} \label{LEM:expectation}
For sufficiently small $\beta$, 
we have 
\begin{align}
& \mathbb{E}_{\mu_\beta} \big[\int :u^2  \!  :_\beta \big] 
 = 0,     \quad
 \hspace{13.5pt} 
\mathbb{E}_{\mu_\beta} \Big[ \big( \int :u^2 \!  :_\beta \big)^2 \Big] \sim \beta^{-\frac{1}{2}},
 \label{expectation1} \\
 & \mathbb{E}_{\mu_\beta} \big[\int :u^4 \!  :_\beta \big]  \lesssim \beta^{-\frac{1}{2}},
  \quad  
\mathbb{E}_{\mu_\beta} \Big[ \big( \int :u^4 \!  :_\beta \big)^2 \Big]    \lesssim \beta^{-\frac{3}{2}} \label{expectation2}.
\end{align}

\end{lemma}

\noi
\begin{proof}[Sketch of Proof] The proof is straightforward for $\int \!  :u^2 \!  :_\beta$.
In view of \eqref{representation1}, we have 
\begin{equation} \label{L4}
 \int u^4 = \sum_{\substack{n_{1234} = 0\\ \ n_j \ne 0}} 
\prod_{j = 1}^4 \frac{g_{n_j}}{\sqrt{ 1 +  \beta n_j^2 }},
\end{equation}

\noi
where $n_{1234}:= n_1 + \cdots + n_4$.
In taking an expectation, the only contribution comes from 
$ n_1 = -n_2$, $n_3 = -n_4$ up to permutations of the indices.
This gives the first term in \eqref{expectation2}.

Now, we say that we have a ``pair'' if we have $n_j = -n_k$, $j \ne k$ in the summation in \eqref{L4}.
If we have a pair, say $n_1 = -n_2$, then we also have $n_3 = -n_4$ since $n_{1234} = 0$.
Then, we can separate the sum in \eqref{L4} as
\begin{align} \label{L4sum}
 \sum_{\substack{n_{1234} = 0\\ \ n_j \ne 0}} 
 = \sum_{\text{pair}} + \sum_{\text{no pair}} 
& = \sum_{\substack{n_1 = -n_2, \ n_3 = -n_4\\ n_j \ne 0 }} + \sum_{\substack{n_1 = -n_3, \ n_2 = -n_4\\ n_j \ne 0 }}
+ \sum_{\substack{n_1 = -n_4, \ n_2 = -n_3\\ n_j \ne 0 }} + \sum_{\text{no pair}} \notag \\
& = 3 \sum_{\substack{n_1 = -n_2, \ n_3 = -n_4\\ n_j \ne 0 }} 
 + \sum_{\text{no pair}} +\text{ error terms}
\end{align}

\noindent
by symmetry. 
The error terms appear from the intersections of the events such as 
$\{n_1 = -n_2\} \cap \{n_1 = -n_3\}$.
They do not have any significant contribution,
and we drop them in the following.\footnote{By precisely computing the error terms, 
one can indeed show that $ \mathbb{E}_{\mu_\beta} \big[\int :u^4 \!  :_\beta \big] =0$.
See \cite{OQV} for details.}
Then, from \eqref{Wick4}, we have
\begin{align} 
\int :  u^4\!  :_\beta  
 \ = 12 \bigg( \sum_{n \geq 1} \frac{|g_{n}|^2 - 1}{1 + \beta n^2}\bigg)^2
 + \sum_{\text{no pair}} \prod_{j = 1}^4 \frac{g_{n_j}}{\sqrt{ 1 +  \beta n_j^2 }} \label{ZWick4}
=:  \I + \II.
\end{align}

\noi
By direct computation, we have 
$\mathbb{E} [\I^2]  \lesssim \beta^{-1}$
and $\mathbb{E}[\I \cdot\II] = 0$.
Finally, we consider
\begin{align*}
\mathbb{E}[\II^2]
=  \mathbb{E} \bigg[\bigg( \sum_{\substack{n_{1234} = 0\\ n_j \ne 0\\ \text{no pair}}} \prod_{j = 1}^4 
\frac{g_{n_j}}{\sqrt{ 1 +  \beta n_j^2 }} \bigg)
\bigg( \sum_{\substack{k_{1234} = 0\\ k_j \ne 0\\ \text{no pair}}} 
\prod_{j = 1}^4 \frac{g_{k_j}}{\sqrt{ 1 +  \beta k_j^2 }} \bigg)\bigg].
\end{align*}

\noi
Note that the indices $\{ n_j\}$ and $\{ k_j\}$ contain no pair
and that $\mathbb{E} [g_n^j] = 0$ for complex-valued Gaussians.
Hence, the  only nonzero contribution comes from $\{ n_1, n_2, n_3, n_4\} = - \{k_1, k_2, k_3, k_4\}$. 
Now, by further separating the summation into (a) $n_j$ all distinct, (b) $n_1 = n_2 \ne n_3, n_4$ and $n_3 \ne n_4$, 
and (c) $n_1 = n_2 = n_3 \ne n_4$ (up to permutations of the indices),
we see that the main contribution comes from  (a) $n_j$ all distinct,
which yields the second term in  \eqref{expectation2}
by a simple Riemann sum approximation.
\end{proof}

We point out that the ``no pair, all distinct"
is responsible for the largest contribution,
which appears again in Subsection 4.2.
Now, define $\mathcal{A}_{\beta, N}$ and $\mathcal{B}_{\beta, N}$ by
\begin{align} 
\mathcal{A}_{\beta, N}   = \big\{  \Big| \int : u^4 \! :_\beta \Big| \leq N \beta^{-\frac{3}{4}}\big\},  
\text{ and }
\mathcal{B}_{\beta, N}  = \big\{  \Big| \int :u^2 \!:_\beta \Big|  \leq N \beta^{-\frac{1}{4}} \big\} \label{setAB}
\end{align}

\noi
for large $N$ and small $\beta > 0$, and we consider the contributions from 
\[ \text{(i) }\mathcal{A}_{\beta, N} \cap \mathcal{B}_{\beta, N}, \quad
\text{ (ii) } \mathcal{A}_{\beta, N} \cap \mathcal{B}^c_{\beta, N}, \quad
\text{ and \quad (iii) } \mathcal{A}^c_{\beta, N}.\]

\noi
First, note that by Chebyshev's inequality with Lemma \ref{LEM:expectation} and \eqref{setAB}, 
we have
\begin{equation}\label{small}
\int_{\mathcal{A}_{\beta, N}^c \cup \mathcal{B}^c_{\beta, N}} d \mu_\beta \lesssim N^{-2}.
\end{equation}

\noi
Hence, we expect that the main contribution for the weak convergence \eqref{mainconv} comes from (i),
and that the contributions from (ii) and (iii) are small.

\smallskip

\noi
$\bullet$ {\bf (i)} On  $\mathcal{A}_{\beta, N} \cap \mathcal{B}_{\beta, N}$:
Now, for sufficiently small $\beta > 0$, 
consider the Gaussian measure $d \wt{\mu}_\beta = \exp \{6 \beta  a_\beta \int  u^2\} \, d \mu_\beta$
(with appropriate normalization.)
Then, under $\wt{\mu}_\beta$, we have
\begin{equation} \label{representation2}
u(x) = \sum_{n \ne 0} \frac{g_n}{\sqrt{1 - 12 \beta  a_\beta + \beta n^2}} e^{2\pi i n x}.
\end{equation}

\noi
Note that $ 12 \beta  a_\beta \sim \beta^\frac{1}{2} \to 0 $ as $\beta \to 0$
and it does not cause a problem.
Then, it converges to the white noise.

\begin{lemma} \label{LEM:conv1}
There exists $C_\beta$ and  $\wt{C}_\beta$ such that 
 we have 
\begin{equation} \label{conv1}
\lim_{\beta \to 0}  C_\beta  \int e^{i \int f u  + 6 \beta a_\beta \int u^2 - 3 \beta a_\beta^2 } d \mu_\beta 
= \lim_{\beta \to 0}  \wt{C}_\beta  \int e^{i \int f u  - 3 \beta a_\beta^2 } d \wt{\mu}_\beta 
=  e^{-\frac{1}{2} \|f\|_{L^2}^2}, 
\end{equation}

\noi
for any smooth mean 0 function $f$ on $\T$,
\end{lemma}

\noi
This follows from a direct computation:
\begin{align*}
\int e^{i \int f u}  d \wt{\mu}_\beta 
 = \exp \Big\{ -\frac{1}{2}  \sum_{n \ne 0} \frac{|\ft{f}_n|^2 }{1 - 12 \beta  a_\beta + \beta n^2}  \Big\} 
\to e^{-\frac{1}{2} \|f \|_{L^2}^2}.
\end{align*}

Next, we show that $\beta \int u^4$ is small in this case 
and that it does not affect the weak convergence in Lemma \ref{LEM:conv1}.
For conciseness of the presentation, 
let \[I_f(F) = \int F(u)  e^{i \int f u  + 6 \beta a_\beta \int u^2 - 3 \beta a_\beta^2 } d \mu_\beta.\]

\begin{lemma} \label{LEM:bound3}
\begin{align} \label{bound3}
\limsup_{\beta \to 0} \bigg|\int_{\mathcal{A}_{\beta, N}\cap \mathcal{B}_{\beta, N}}  
\chi_{\{ \int u^2 \leq K \beta^{-\frac{1}{2}}\}} & e^{i \int f u +  \beta \int u^4} d \mu_\beta 
- I_f(1)\bigg| 
\lesssim  N^{-1} .
\end{align}
\end{lemma}

\begin{proof}
On $\mathcal{A}_{\beta, N}$, we have 
$ \big|e^{\beta \int : u^4  :_\beta } - 1\big|\lesssim \beta^\frac{1}{4} N $ for $\beta \leq  N^{-4}$.
Hence, we have 
\begin{align*}
 \bigg|\int_{\mathcal{A}_{\beta, N}\cap \mathcal{B}_{\beta, N}}  
 & \chi_{\{ \int u^2 \leq K \beta^{-\frac{1}{2}}\}}  e^{i \int f u +  \beta \int u^4} d \mu_\beta 
  - I_f\big( \chi_{\mathcal{A}_{\beta, N}\cap \mathcal{B}_{\beta, N}}
  \chi_{ \{ \int u^2 \leq K \beta^{-\frac{1}{2}}\}}\big) \bigg| \\
&  \lesssim    e^{ 6 \beta^\frac{1}{2} a_\beta K - 3\beta a_\beta^2}   \int | e^{\beta \int : u^4  :_\beta } - 1| d \mu_\beta
\lesssim  \beta^\frac{1}{4} N.
\end{align*}

\noi
Note that we have $\mathcal{B}_{\beta, N} \subset \{ \int u^2 \leq K \beta^{-\frac{1}{2}}\}$ for sufficiently small $\beta$.
Hence,  it suffices to show
\begin{equation}\label{bound31}
\limsup_{\beta \to 0}
|I_f( \chi_{\mathcal{A}_{\beta, N}\cap \mathcal{B}_{\beta, N}})
-I_f(1)| 
=
\limsup_{\beta \to 0}
|I_f( \chi_{\mathcal{A}^c_{\beta, N}\cup \mathcal{B}^c_{\beta, N}})|
\lesssim N^{-1}.
\end{equation}

\noi
By Cauchy-Schwarz inequality along with \eqref{small}, we have
\begin{align} \label{ABC}
|I_f( \chi_{\mathcal{A}^c_{\beta, N}\cup \mathcal{B}^c_{\beta, N}})|
\leq \mu_\beta^\frac{1}{2} (\mathcal{A}^c_{\beta, N}\cup \mathcal{B}^c_{\beta, N})
\bigg(\int  e^{   6 \beta a_\beta \int u^2  } d \mu_\beta\bigg)^\frac{1}{2}
\lesssim N^{-1}
\end{align}

\noi
since 
$e^{   6 \beta a_\beta \int u^2  } d \mu_\beta$ is a normalizable 
density for $\beta$ small in view of  $\beta a_\beta\to 0$ as $\beta \to 0$.
\end{proof}

\smallskip
\noi
$\bullet$ {\bf (ii)} On $\mathcal{A}_{\beta, N} \cap \mathcal{B}^c_{\beta, N}$:
On $ \mathcal{A}_{\beta, N}\cap\{ \int u^2 \leq K \beta^{-\frac{1}{2}}\}$,
we have
\[ \beta \int u^4 \leq \beta \Big|\int :u^4 \!:_\beta\Big| + 6 \beta a_\beta \int u^2 + 3 \beta a_\beta^2 
\lesssim 1 \]

\noi
for $\beta \leq N^{-4}$.
Hence, by \eqref{small}, we have
\begin{equation} \label{bound2}
\int_{\mathcal{A}_{\beta, N}\cap \mathcal{B}^c_{\beta, N}}  
\chi_{\{ \int u^2 \leq K \beta^{-\frac{1}{2}}\}} e^{i \int f u +  \beta \int u^4} d \mu_\beta \lesssim N^{-2} .
\end{equation}

\smallskip
\noi
$\bullet$ {\bf (iii)} On $\mathcal{A}^c_{\beta, N}$:
In this case, we have
\begin{equation*} 
\int_{\mathcal{A}^c_{\beta, N}}  \chi_{\{ \int u^2 \leq K \beta^{-\frac{1}{2}}\}} 
e^{i \int f u +  \beta \int u^4} d \mu_\beta \lesssim N^{-1} .
\end{equation*}
This follows from Cauchy-Schwarz inequality followed by \eqref{ExpEx2} and \eqref{small}.
Hence, it remains to prove the exponential expectation \eqref{ExpEx2},
which is by far the most technical part of the proof.
Note that the general framework for the $p = 3 $ case is similar, but simpler.

\subsection{Exponential expectation}

It suffices to show the tail estimate
\begin{equation}\label{Tail1}
\mu_\beta \big[ \, \beta \int u^p > \ld, \int u^2 \leq K \beta^{-\frac{1}{2}} \big] 
\leq e^{-c \ld^{1 + \dl}}
\end{equation}

\noi
for $p = 4$, uniformly in small $\beta > 0$, 
where $\mu_\beta$ is as in \eqref{Gauss2}.
For $\beta = 1$, Bourgain \cite{BO4} proved
\eqref{Tail1} for $p<6$
via the dyadic pigeonhole principle with the large deviation lemma (Lemma \ref{LEM:polar})
as in the proof of Proposition \ref{PROP:meas}.
We point out that Bourgain's argument is not sufficient even for $p=3$. See \cite{OQV}.

Following Bourgain's argument, 
we can prove \eqref{Tail1}
for 
\begin{itemize}
\item all $\ld > 0$ on $\{ |n| \geq \beta^{-1-}\}$, i.e. large frequencies
\item $\ld \geq \beta^{-\frac{1}{2}-}$ with no frequency restriction.
\end{itemize}

\noi
Hence, we need to show \eqref{Tail1} for $\ld \leq \beta^{-\frac{1}{2}-}$,
assuming that $u$ has a {\it finite} Fourier support.

First, note that we have
$\beta \int u^4 = \beta \int \! : \! u^4 \!:_\beta + O(1)$ on $\{ \int u^2 \leq K \beta^{-\frac{1}{2}}\}$.
We prove \eqref{Tail1} with $\beta \int \! : \! u^4 \!:_\beta$ instead of $\beta \int u^4$.  
As before, the main contribution comes from ``no pair, all distinct".
In the following, we prove
\begin{equation}\label{Tail2}
\mu_\beta \big[ |Q_\beta| > \ld, \int u^2 \leq K \beta^{-\frac{1}{2}} \big] 
\leq e^{-c \ld^{1 + \dl}}
\end{equation}

\noi
for $\ld \leq \beta^{-\frac{1}{2}-}$,  where
\begin{equation} \label{QQ}
Q_\beta = \beta \sum_{**} \prod_{j = 1}^4 \frac{g_{n_j}}{\sqrt{ 1 +  \beta n_j^2 }} 
\end{equation}
with
$** = \{ \, n_{1234} := n_1 + \cdots + n_4 = 0,   \text{ no pair, all distinct, } |n_j|\leq \beta^{-1-}\}$.

\smallskip

Now, we give a brief review on the hypercontractivity of the Ornstein-Uhlenbeck semigroup.
See Tzvetkov \cite[Sec.3]{TZ3} for details.
Let $L$ denote the generator of the Ornstein-Uhlenbeck process on $H := L^2(\R^d, e^{-|x|^2/2}dx)$
given by $L = \Dl - x \cdot \nabla$.
Then, let $S(t) = \exp( tL)$
be the semigroup associated with $\dt u = L u$.
Then, the hypercontractivity of the Ornstein-Uhlenbeck semigroup
says the following:
\begin{lemma} \label{LEM:HYP1}
Let $q \geq 2$. For $f \in H$ and $t\geq \frac{1}{2} \log(q-1)$,
we have
\[ \| S(t) f \|_{L^q(\R^d, \exp(-|x|^2/2)dx)}
\leq \|  f \|_{L^2(\R^d, \exp(-|x|^2/2)dx)}\]
\end{lemma}

The eigenfunctions of $L$ are given by $\prod_{j = 1}^d h_{k_j}(x_j)$, 
where $h_k$ is the Hermite polynomial of degree $k$,
and the corresponding eigenvalue is given by
 $\ld = - (k_1+ \cdots + k_d)$.
We list the first few Hermite polynomials:
\begin{equation} \label{Hermite}
h_0(x) = 1,  \ h_1(x) = -x, \ h_2(x) = \tfrac{1}{\sqrt{2}}(x^2 - 1), \ \cdots
\end{equation}

\noi
For our application, 
let \[H(x) = \sum_\G c(n_1, \cdots, n_4) x_{n_1} \cdots x_{n_4},\] 

\noi
where 
$\G = \{ (n_1, \cdots, n_4) \in \{1, \cdots, d\}^4, \text{ all distinct} \}.$
Note that $H(x)$ is an eigenfunction of $L$ with the eigenvalue $-4$.
Then, by Lemma \ref{LEM:HYP1}, we have
the following {\it dimension-independent} estimate:
\begin{equation} \label{LEM:HYP2}
 \| H(x) \|_{L^{q}(\R^d, \exp(-|x|^2/2)dx)}
\leq q^2 \|  H(x)  \|_{L^2(\R^d, \exp(-|x|^2/2)dx)}.
\end{equation}

By expanding the complex-valued Gaussians $g_n$ into their real and imaginary parts, 
we can apply \eqref{LEM:HYP2} to $Q_\beta$ in \eqref{QQ}.
From (the proof of) Lemma \ref{LEM:expectation}, we have
$\|Q_\beta\|_{L^2(d\mu_\beta)} \leq C \beta^{\frac{1}{4}}$.
By \eqref{LEM:HYP2}, 
we have \begin{equation} \label{QQQ}
\|Q_\beta\|_{L^q(d\mu_\beta)} 
\leq C q^2 \beta^\frac{1}{4}
\end{equation}

\noi
for all $ q\geq 2$.
It is important that $u$ has a finite Fourier support, but 
the actual upperbound on the support is not important.
Then, we have
\begin{equation}\label{QQQQ}
 \int \exp( c \beta^{-\frac{1}{8}} |Q_\beta|^\frac{1}{2} ) d\mu_\beta \leq C
 \end{equation}

\noi 
from Lemma 4.5 in \cite{TZ3}
(or equivalently, expanding the exponential in the Taylor series
and applying \eqref{QQQ}.)
\eqref{QQQQ} in turn implies
$ \mu_\beta [ |Q_\beta| > \ld ] \leq \exp( -c' \beta^{-\frac{1}{8}} \ld^\frac{1}{2} )$, 
i.e.
we proved \eqref{Tail2} for $\ld \leq \beta^{-\frac{1}{4}+}$.

Now, we consider the remaining case: $ \beta^{-\frac{1}{4}+} \leq \ld \leq \beta^{-\frac{1}{2}-}$.
Then, using $\ld \geq \beta^{-\frac{1}{4}+\eps}$,
\begin{align*}
\mu_\beta \big( |Q_\beta| \geq \ld ) & \leq \frac{\|Q_\beta\|^q_{L^q(d \mu_\beta)}}{\ld^q}
\leq C q^{2q} \beta^{\frac{q}{2} -\eps q}
\leq e^{2q \ln q} e^{-\frac{q}{3} \ln \beta^{-1}}
= e^{ - \frac{q}{3} \ln \beta^{-1} + 2q \ln q}
\intertext{By choosing $q \sim \beta^{-\frac{3}{4}} \ll \beta^{-1}$ and using $\ld \leq  \beta^{-\frac{1}{2}-\eps}$, }
& \leq e^{ - c \beta^{-\frac{3}{4}} \ln \beta^{-1}}
\leq e^{-c\ld^{\frac{3}{2}-}}.
\end{align*}

\noi
This completes the proof of the tail estimate \eqref{Tail1},
and hence the exponential expectation \eqref{ExpEx2}.

\section{Remarks and Comments}

\noi
(a) We summarize the different approaches we discussed in this paper.

\smallskip
\noi
$\bullet$
Method 0,  Complete integrability approach: 
It uses strong results
which are only true for KdV, and thus it can not be applied to non-integrable KdV variants.

\smallskip
\noi
$\bullet$
Method 1, PDE approach: 
It is a direct approach, only requiring local well-posedness on the support 
of the white noise.
It can be applied to non-integrable KdV variants as well.
However, one needs to establish  LWP with a good estimate
which is often nontrivial.

\smallskip
\noi
$\bullet$
Method 2, Probabilistic approach:
This can be used to establish a formal invariance
even when  well-posedness is not known.
It can also establish the invariance as soon as there is a continuous flow for the PDE.
i.e. it does not require any PDE estimate.
However, one needs to have a continuous flow which needs to be proven elsewhere.

\medskip
\noi
(b) We established the formal invariance of the white noise for mKdV and 1-$d$ cubic NLS.
However, the existence of a continuous flow in the support of white noise,
 which is needed to prove the invariance,
is not known.
Note that it seems essential to study the 
{\it Wick ordered} cubic NLS below $L^2(\T)$
in place of the usual cubic NLS.
See Remark \ref{RM:2.3}.
In this respect, mKdV seems much harder to treat.
Takaoka-Tsutsumi \cite{TT} proved LWP of the Wick ordered mKdV
in $H^s(\T)$ for $s > \frac{3}{8}$.
However, this is far from reaching the support of the white noise.

There are several partial results of the well-posedness of the  1-$d$ Wick ordered cubic NLS  outside $L^2(\T)$.
Christ \cite{Christ} constructed local-in-time solutions  in $\mathcal{F}L^{0, p}$ for $p<\infty$
by his power series method.
Also see Gr\"unrock-Herr \cite{GH} for the same result via the fixed point argument.
Colliander-Oh \cite{CO1} constructed local-in-time solutions
with Gaussian-randomized initial data on the negative Sobolev spaces 
by exhibiting nonlinear smoothing under randomization as in \cite{BO7}.
The proof is probabilistic and uses the estimates on the homogeneous Wiener chaos
as a result of the hypercontractivity of the Ornstein-Uhlenbeck semigroup.
Then, such local-in-time solutions are extended to global ones
(in the absence of invariant measures)
via the so-called Bourgain's high-low method \cite{BO98}.

In \cite{BO7} and \cite{CO1}, the local solutions were
constructed via the fixed point argument around the linear solution
$z_1(t) := S(t)u_0$
with probabilistic arguments.
Also see Burq-Tzvetkov \cite{BT2} and  Thomann \cite{LT} for related arguments.
While the basic probabilistic argument
is similar, the argument in \cite{BT2, LT} further exploits the properties of the eigenfunctions,
and the argument in \cite{BO7, CO1} exploits more properties of the product of Gaussians
via the hypercontractivity of the Ornstein-Uhlenbeck semigroup.
%
Such construction for KdV fails on the support of the white noise.
Nonetheless, in Oh \cite{OH7},
local-in-time solutions are constructed 
via the nonlinear analysis on the second iteration.
See Section 3 and \cite{BO3, OH5, OH6}.

Another possible improvement is to construct
solutions via the fixed point argument 
around the second iterate $z_2(t) := S(t) u_0 + c \int_0^t S(t - t') \mathcal{N}(S(t)u_0) dt'$,
where $\mathcal{N}(u)$ is the nonlinearity of a PDE.
This idea may be useful to study the Gibbs measure for the Benjamin-Ono equation.
On the one hand, Tzvetkov \cite{TZ3} constructed the Gibbs measure for the Benjamin-Ono equation,
which is supported in $\bigcap_{s <0} H^s (\T) \setminus L^2(\T)$.
On the other hand, Molinet \cite{MO1, MO2} proved the sharp well-posedness in $L^2(\T)$.
It was also shown in \cite[Proposition 5.4]{TZ3}
that $z_2(t) - z_1(t)$ is not in $L^2(\T)$.
While one can not construct solutions around the linear solution
(since $z_2(t) - z_1(t) \notin L^2(\T)$),
it seems reasonable to try to construct solutions on the support of the Gibbs measure
via the fixed point argument around the second iterate $z_2(t)$.
\smallskip

\noi
{\bf Acknowledgments:} 
The author would like to thank Prof. Henry P. McKean 
for telling him about this subject.
It has taken him more than several years to digest, but 
different results emerged from 
 the conversation we had at NYU.

\end{document}